\documentclass[10pt]{article}

\usepackage{bm}
\usepackage{fullpage}  

\usepackage{url}
\usepackage{algpseudocode}
\usepackage{multirow} 
\usepackage{hhline}  
\usepackage{verbatim}
\usepackage{graphicx}
\usepackage{amsmath,amssymb,amsfonts,graphicx,amsthm,mathtools,nicefrac}
\usepackage{lscape}
\usepackage{color}
\usepackage{authblk}
\usepackage{subfigure}
\usepackage[ruled,vlined]{algorithm2e}

\usepackage{bbding}
\usepackage{indentfirst}
\usepackage{multirow}

\allowdisplaybreaks
\newcommand\numberthis{\addtocounter{equation}{1}\tag{\theequation}}

\newtheorem{theorem}{Theorem}[section]
\newtheorem{assumption}{Assumption}[section]
\newtheorem{remark}{Remark}[section]

\numberwithin{equation}{section}

\newtheorem{lemma}[theorem]{Lemma}

\def\la {\left\langle}
\def\ra {\right\rangle} 

\newcommand{\matsnorm}[2]{\left\| #1\right\|_{{#2}}}

\newcommand{\twonorm}[1]{\ensuremath{\matsnorm{#1}{\footnotesize{2}}}}

\newcommand{\bfm}[1]{\bm{#1}}

\newcommand{\E}[2][]{\mathbb{E}_{#1} \left\{ #2 \rule{0mm}{3mm}\right\}}
 
\def\vg{\bfm g}

     \def\R{\mathbb{R}}
     \def\S{\mathbb{S}}

\def \vtheta {\bfm \theta}

\def\calA{{\cal  A}}

\def\calS{{\cal  S}}

\def \tran {\mathsf{T}}

\usepackage{bbm}

\newcommand{\kw}[1]{{#1}}
\begin{document}
	
\title{On the Linear Convergence of Policy Gradient under Hadamard Parameterization}
\author[1]{Jiacai Liu}
\author[2]{Jinchi Chen}
\author[1]{Ke Wei}
\affil[1]{School of Data Science, Fudan University, Shanghai, China.}
\affil[2]{School of Mathematics, East China University of Science and Technology, Shanghai, China.\vspace{.15cm}}
\date{\today}

\maketitle
\begin{abstract}
The convergence of deterministic policy gradient  under the Hadamard parameterization is studied in the tabular setting and the linear convergence of the algorithm is established. To this end, we first show that the error decreases at an $O(\frac{1}{k})$ rate for all the iterations. Based on this result, we further show that the algorithm has a faster local linear convergence rate after $k_0$ iterations, where $k_0$ is a constant that only depends on the MDP problem and the  initialization. 
\kw{To show the local linear convergence of the algorithm, we have indeed established the  contraction of the sub-optimal probability  $b_s^k$ (i.e., the probability of the output policy $\pi^k$ on non-optimal actions)  when $k\ge k_0$.}\\

\end{abstract}

\section{Introduction}
Reinforcement learning (RL) is a type of machine learning technique for solving sequential decision problems and has recently achieved great success in many areas, such as games \cite{mnih2015humanlevel,SilverHuangEtAl16nature,Berner2019Dota2W,StarCraft}, robotics \cite{robot1,robot2,robot3} and various other real applications \cite{Agarwal2016MakingCD,chen-2019-topk,chipdesign}. Typically, RL can be modeled as a Markov decision process (MDP) represented by a tuple $\mathcal{M} \left(\calS ,\calA ,P,r,\gamma ,\mu \right)$, where $\calS$ is the state space, $\calA$ denotes the action space, $P(s'|s,a)$ is the transition probability or density from state $s$ to state $s'$ under action $a$, $r: \calS \times \calA \times \calS \rightarrow \mathbb{R}$ is the reward function, $\gamma \in [0,1)$ is the discounted factor  and $\mu$ is the probability distribution of the initial state $s_0$. Here we restrict our attention to the tabular setting where $\calS$ and $\calA$ are finite, i.e.  $|\calS|<\infty$ and $|\calA|<\infty$. 
Let $\Delta(\calA)$ be the probability simplex over the set $\calA$, 
defined as
\begin{align}
\Delta \left( \calA \right) =\left\{ \theta \in \mathbb{R} ^{|\calA |}: \theta _i\ge 0,\sum_{i=1}^{|\calA |}{\theta _i}=1 \right\}.
\end{align}
Given a stationary policy $\pi:\calS\rightarrow \Delta(\calA)$, the state value function at  $s\in\calS$ and the state-action value function at $(s,a)\in\calS\times \calA$ are defined as
\begin{align}
V^{\pi}\left( s \right) &:=\mathbb{E} \left\{ \sum_{t=0}^{\infty}{\gamma ^tr\left( s_t,a_t,s_{t+1} \right)}|s_0=s,\pi \right\},
\label{value function} \\
Q^{\pi}\left( s,a \right) &:=\mathbb{E} \left\{ \sum_{t=0}^{\infty}{\gamma ^tr\left( s_t,a_t,s_{t+1} \right)}|s_0=s,a_0=a,\pi \right\}.
\label{Q function}
\end{align}
Overall, the goal of RL is to a find a policy that maximizes the weighted average of the state values under the initial distribution $\mu$, namely to solve
\begin{align}
\underset{\pi}{\max}\,V^{\pi}\left( \mu \right). \label{objective}
\end{align}
Here $
V^{\pi}\left( \rho \right) :=\mathbb{E} _{s\sim \rho}\left\{ V^{\pi}\left( s \right) \right\}$ for any $\rho \in \Delta(\calS)$.

Overall, there are two typical categories of RL algorithms. Value-based methods, such as value iteration \cite{puterman2014markov}, Q-learning \cite{Qlearning} and DQN \cite{mnih2015humanlevel}, seek the optimal policy based on the idea of fixed point iteration or approximate dynamic programming. In contrast, policy-based methods, including policy gradient \cite{puterman2014markov}, TRPO \cite{trpo} and PPO \cite{ppo} reformulate \eqref{objective} as an optimization problem based on  differentiable parameterized policies and then solve the problem by optimization methods.

In this work, we focus on the policy optimization based methods \cite{pg1,pg2,pg3}. 
Using a differentiable parameterized  policy $\pi_\theta$ where $\theta \in \R^d $, \eqref{objective} can be written as a finite dimensional optimization problem:
\begin{align}
\underset{\theta \in \mathbb{R} ^d}{\max}\,V^{\pi _{\theta}}\left( \mu \right) .
\label{PG Objective}
\end{align}
The most widely studied parameterizations in tabular setting are direct parameterization:
\begin{align}
\pi _{\theta}\left( a|s \right) =\theta _{s,a} \quad \text{s.t.}\quad \sum_{a\in \calA}{\theta _{s,a}}=1,
\label{simplex}
\end{align}
and the softmax parameterization:
\begin{align}
\pi _{\theta}\left( a|s \right) =\frac{\exp \left\{ \theta_{s,a} \right\}}{\sum_{a\in \calA}{\exp \left\{ \theta_{s,a} \right\}}}.
\label{softmax}
\end{align}
Policy gradient (PG) is a first order method for solving \eqref{PG Objective}.
Besides, natural policy gradient (NPG) \cite{NIPS2001_4b86abe4} is an important variant  of PG,  which can be viewed as a preconditioned first order method \cite{10.1162/089976698300017746}. Moreover, it is shown that the NPG method can be cast as a special case of the policy mirror dscent (PMD) method \cite{DBLP:conf/aaai/ShaniEM20}.

\subsection{Motivation}
A recent work \cite{Li_McKenzie_Yin} in optimization  shows that converting an optimization problem on a standard simplex into one on a unit sphere may have several advantages, such as \kw{the projection onto the unit sphere admits a simple closed-form solution} and the corresponding algorithm can be significantly faster  empirically. In the RL setting, the objective \eqref{objective} can be seen as an optimization problem on the Cartesian product of simplices. Inspired by \cite{Li_McKenzie_Yin}, we consider the Hadamard parameterization \kw{for the RL problem}:
\begin{align}
\label{Hadamard Policy}
    \pi_{\theta}(\cdot|s) = \theta_s \odot \theta_s\quad \text{s.t.} \quad \theta_s\in\S^{|\calA|-1},
\end{align}
\kw{which leads to an optimization problem} constrained on the Cartesian product of unit spheres. Hadamard PG utilizes the Riemannian gradient ascent to solve this problem. As far as we know,  \kw{the study of the Hadamard parameterization has not been well-received in RL.} A most related work is \cite{Mei_Xiao_Dai_Li_Szepesvari_Schuurmans_2020}, which studies the escort policy
\begin{align}
\label{Escort Policy}
    \pi_{\theta}(a|s) = \frac{\left|\theta_{s,a}\right|^p}{\sum_{a^\prime \in \calA}{\left|\theta_{s,a^\prime}\right|^p}},
\end{align}
in which the $p=2$ case  can be seen as the normalized version of the Hadamard parameterization.   It is worth noting that the policy sequence $\pi^k$ 
generated by the Hadamard PG with a constant step size $\eta$ is the same as  that generated by the escort PG with $p=2$ and the step sizes $\eta_k(s) = \eta \cdot \left\|\theta_s^k \right\|^2_p$  (\kw{algorithm investigated in \cite{Mei_Xiao_Dai_Li_Szepesvari_Schuurmans_2020};} see Section~\ref{equivalent_algorithm} for more details of the equivalence), where $\theta^k$ is the output parameter at the $k$-th iteration. Thus the sublinear convergence result established in \cite{Mei_Xiao_Dai_Li_Szepesvari_Schuurmans_2020} is also applicable to Hadamard PG. The goal of this paper is  to provide an improved sublinear convergence analysis for Hadamard PG compared with the one in \cite{Mei_Xiao_Dai_Li_Szepesvari_Schuurmans_2020} and meanwhile show that Hadamard PG also enjoys local linear convergence based on the sub-optimal probability analysis framework.

\subsection{Main contributions and related work}
The main contributions of this paper are the improved sublinear convergence of Hadamard PG and  the establishment of its local linear convergence.
 For the sublinear convergence, the result in \cite{Mei_Xiao_Dai_Li_Szepesvari_Schuurmans_2020} implies  that Hadamard PG converges to the optimal value at a rate of  $\mathcal{O}(\frac{\left|\calA\right| \left|\calS\right|}{k(1-\gamma)^6})$. \kw{In contrast}, the rate we have established in this paper is $\mathcal{O}(\frac{1}{k(1-\gamma)^4})$, which has a better dependency on the horizon $\frac{1}{1-\gamma}$ and has no explicit dependency  on the state space size $\left|\calS\right|$ and the action space size $\left|\calA\right|$. Moreover, the local linear convergence is new for Hadamard PG, to the best of our knowledge. The comparisons  with the results in \cite{Mei_Xiao_Dai_Li_Szepesvari_Schuurmans_2020} are summarized in Table \ref{table1}. 

\begin{table}[ht!]
\centering
\caption{Comparison with \cite{Mei_Xiao_Dai_Li_Szepesvari_Schuurmans_2020} on the convergence of Hadamard PG. In the column ``sublinear convergence rate'', $\lambda$ is a positive constant (see Lemma \ref{lemLambda}) and $\tilde{\mu} := \underset{s \in \calS}{\min} \,\,\mu(s)$.}
\vspace{0.1cm}
\label{table1}
\renewcommand{\arraystretch}{2.2}
\begin{tabular}{c|c|c|c} 
\hline
     & Sublinear convergence rate & Linear & Step size  \\
\hline
\cite{Mei_Xiao_Dai_Li_Szepesvari_Schuurmans_2020}  & $\frac{20\left|\calA\right| \left|\calS\right|}{k\cdot(1-\gamma)^6 \cdot \lambda} \cdot \frac{1}{\tilde{\mu}} \cdot \left\|\frac{d^*_\mu}{\mu}\right\|^2_\infty$  & \XSolidBrush  &  $\eta \leq \frac{(1-\gamma)^3}{40\left|\calA\right|}$  \\ 
\hline
\textbf{this paper} & $\frac{4+\kappa^2}{k\cdot 2\kappa(1-\gamma)^4\cdot \lambda} \cdot \frac{1}{\tilde{\mu}^2}$ & \CheckmarkBold  &  $\eta = \frac{(1-\gamma)^2 \kappa}{4}$, $\forall \kappa \in (0,2)$\\
\hline
\end{tabular}
\end{table}

Since the convergence analyses of different policy optimization methods (including exact PG, NPG, PMD) under the simplex and softmax parameterizations have received many investigations recently, we give a short discussion towards this line of research below, with related results being summarized in Table \ref{table2}.

\paragraph{Convergence of PG}
For PG with the constant step size, the best known convergence rate is $O(\frac{1}{k})$ for both the simplex \cite{Agarwal_Kakade_Lee_Mahajan_2019,NEURIPS2020_30ee748d,linxiao2022,LiuLiWei2023PPG} and softmax parameterizations \cite{mei2020global}. Moreover, it has been further shown in \cite{LiuLiWei2023PPG} that \kw{PG enjoys a finite iteration convergence under the simplex parameterization. In contrast,}
 the convergence rate of PG  suffers from an $\Omega \left( \frac{1}{k} \right)$ lower bound for the softmax parameterization.
With exact line search, the authors in \cite{pmlr-v130-bhandari21a} show that PG under the simplex parametrization has a linear convergence rate.  
Different methods have been proposed to improve the convergence of PG under the softmax parameterization. 
 For instance, by introducing entropy regularization into the objective function, the authors in \cite{mei2020global}  show that PG with a constant step size can have a linear convergence rate.   
 Apart from entropy regularization, it is shown in \cite{pmlr-v139-mei21a}  that PG with adaptive step sizes based on the geometric information can also achieve linear convergence.

\paragraph{Convergence of softmax NPG and PMD} 
 Recently, a series of works show that  softmax NPG  enjoys a local or global linear convergence rate \cite{Khodadadian, pmlr-v130-bhandari21a, NPG_chen2021,Zhan_Cen_Huang_Chen_Lee_Chi_2021,linxiao2022,lan2021}. In particular, it is shown in \cite{pmlr-v130-bhandari21a,Khodadadian,Johnson2023}  that softmax NPG  can achieve local linear convergence  with a constant step size, and  global linear convergence with exact line search or adaptive step sizes. Moreover, the authors in \cite{NPG_chen2021,Zhan_Cen_Huang_Chen_Lee_Chi_2021} show that NPG with the entropy regularization has a linear convergence rate globally and  a super-linear convergence rate locally. \kw{As mentioned previously}, softmax NPG can be seen as a special case of PMD (\cite{DBLP:conf/aaai/ShaniEM20})) in the policy space, where the Bregman divergence is the KL divergence.
 In \cite{linxiao2022,Alfano_Yuan_Rebeschini_2023}, it is shown that  PMD enjoys the $O(\frac{1}{k})$ sublinear convergence rate with a constant step size and the  linear convergence rate with  non-adaptive geometrically increasing step sizes.
 For  PMD with a general convex regularizer,    the $O(\frac{1}{k})$ sublinear convergence rate  is established in  \cite{lan2021} under a constant step size. When the regularizer  is strongly convex, \cite{lan2021} also shows that PMD enjoys a linear convergence if the constant step size is large enough. In \cite{lan2021} and \cite{Li_Zhao_Lan_2022}, the global linear convergence of PMD \kw{using a homotopic technique} has been established for the strongly convex regularizer and   non-adaptive geometrically increasing step sizes.  \kw{The convergence of PMD with a wider range of regularizers (including convex and non-smooth regularizers) is discussed in \cite{Zhan_Cen_Huang_Chen_Lee_Chi_2021} and the linear convergence rate is established under a  constant step size.}

\begin{table}[ht!]
\centering
\caption{Summary of current results on the convergence of PG, softmax NPG and PMD. In the table, and ``global'' means whether the algorithm converges to the optimal values of the original MDP.}
\vspace{0.3cm}
\label{table2}
\renewcommand{\arraystretch}{1.5}{\small
\begin{tabular}{c|c|c|c|c|c} 
\hline
Algorithm     & Convergence rate & Global & Linear & Regularizer & Step sizes  \\ 
\hline
Simplex PG   & \multirow{ 2}{*}{$O(\frac{1}{k})$} & \multirow{ 2}{*}{\CheckmarkBold}  & \multirow{ 2}{*}{\XSolidBrush}& \multirow{ 2}{*}{\XSolidBrush} &  \multirow{ 2}{*}{constant}  \\ 
(\cite{NEURIPS2020_30ee748d,Agarwal_Kakade_Lee_Mahajan_2019,linxiao2022,LiuLiWei2023PPG})& & & & &\\
\hline
Softmax PG  (\cite{mei2020global}) & $O(\frac{1}{k})$ & \CheckmarkBold  & \XSolidBrush& \XSolidBrush  &  constant       \\ 
\hline
Softmax PG (\cite{pmlr-v139-mei21a}) &   $O(e^{-c\cdot k})$ & \CheckmarkBold  & \CheckmarkBold& \XSolidBrush  &  adaptive       \\ 
\hline
Softmax PG (\cite{mei2020global})       &  $O(e^{-c\cdot k})$ &\XSolidBrush &    \CheckmarkBold   &    entropy      &  constant\\ 
\hline
\textbf{Hadamard PG}     &\multirow{ 2}{*}{$O((1-\rho(\delta,\kappa,\gamma))^k)$} & \multirow{ 2}{*}{\CheckmarkBold}  & \multirow{ 2}{*}{\CheckmarkBold} &  \multirow{ 2}{*}{\XSolidBrush}  &  \multirow{ 2}{*}{constant}  \\
(\textbf{this paper})& & & & & \\
\hline
Escort PG \cite{Mei_Xiao_Dai_Li_Szepesvari_Schuurmans_2020} & $O(\frac{1}{k})$ & \CheckmarkBold  & \XSolidBrush &  \XSolidBrush  &  constant  
\\ 
\hline
Softmax NPG  (\cite{Khodadadian})  & $O\left( e^{-k\left( 1-\frac{1}{\lambda}\eta \right) \Delta} \right)$ & \CheckmarkBold   &    \CheckmarkBold          &      \XSolidBrush   &   constant \\ 
\hline
Softmax NPG (\cite{pmlr-v130-bhandari21a})  & $O\left( \left( \frac{1+\gamma}{2} \right) ^k \right)$ & \CheckmarkBold   &    \CheckmarkBold          &      \XSolidBrush   &   adaptive \\ 
\hline
Softmax NPG (\cite{Johnson2023})  & $O(\gamma^k)$ & \CheckmarkBold   &    \CheckmarkBold          &      \XSolidBrush   &   adaptive \\ 
\hline
Softmax NPG  (\cite{NPG_chen2021})      &   $O\left(\left( 1-\eta \tau \right) ^k \right)$  &\XSolidBrush &   \CheckmarkBold     &      entropy     &    constant       \\
\hline
PMD (\cite{linxiao2022})    &   $O\left(\left( 1 - \frac{1}{\vartheta_\rho} \right) ^k \right)$  &\CheckmarkBold &   \CheckmarkBold     &       \XSolidBrush      &    non-adaptive increasing  \\
\hline
PMD (\cite{Alfano_Yuan_Rebeschini_2023})    &   $O((1-\frac{1}{C_3})^k)$  &\CheckmarkBold &   \CheckmarkBold     &       \XSolidBrush      &    non-adaptive increasing  \\
\hline
PMD (\cite{lan2021})    &   $O(\gamma^k)$  &\XSolidBrush &   \CheckmarkBold     &      strongly convex      &    $1+\eta\mu \ge \frac{1}{\gamma}$\\
\hline
APMD \cite{lan2021}     &   \multirow{ 2}{*}{$O(\gamma^k)$}  &\multirow{ 2}{*}{\CheckmarkBold} &   \multirow{ 2}{*}{\CheckmarkBold}     &     \multirow{ 2}{*}{diminishing} &   \multirow{ 2}{*}{non-adaptive increasing}\\
\& HPMD \cite{Li_Zhao_Lan_2022}&&&&&\\
\hline
PMD (\cite{Zhan_Cen_Huang_Chen_Lee_Chi_2021})    &   $O((1-(1-\frac{1}{\eta \tau})(1-\gamma))^k)$  &\XSolidBrush &   \CheckmarkBold     &   convex &   constant\\
\hline
\end{tabular}}
\end{table}
\subsection{More Notations and Organization}\label{sec:notation}

\paragraph{Advantage function.}
Recalling the definitions of the state value function (\ref{value function}) and the state-action value function (\ref{Q function}), the advantage function of a policy $\pi$ is defined as
\begin{align*}
A^{\pi}\left( s,a \right) :=Q^{\pi}\left( s,a \right) -V^{\pi}\left( s \right).
\end{align*}
The advantage function $A^\pi(s,a)$ measures how well a single action is  compared with the current  policy $\pi$.

\paragraph{State visitation probability.}  The discounted state visitation probability from any initial state $s_0 \in \calS$ is defined as 
\begin{align*}
d_{s_0}^{\pi}\left( s \right) :=\left( 1-\gamma \right) \mathbb{E} \left\{ \sum_{t=0}^{\infty}{\gamma ^t\mathbf{1} [ s_t=s ]}|s_0,\pi \right\},
\end{align*}
where $\mathbf{1}[\cdot]$ is the indicator function. Moreover, define $d_{\rho}^{\pi}\left( s \right) :=\mathbb{E} _{s_0\sim \rho}\left\{ d_{s_0}^{\pi}\left( s \right) \right\}$ for any $\rho \in \Delta(\mathcal{S})$. A direct computation yields that
\begin{align}
d_{\rho}^{\pi}\left( s \right) \ge \left( 1-\gamma \right) \mathbb{E} _{s_0\sim \rho}\{\mathbf{1}[s_0=s] \}\ge \left( 1-\gamma \right) \min_{s\in \mathcal{S}} \,\rho (s)\quad \mbox{for any policy } \pi.
\label{d_eq}
\end{align}

\paragraph{Optimal values and optimal policies.} The optimal values are defined as the supremum over the all possible polices, i.e. 
\begin{align}
V^*\left( s \right) :=\underset{\pi}{\mathrm{sup}}\,V^{\pi}\left( s \right),\quad Q^*\left( s,a \right) :=\underset{\pi}{\mathrm{sup}}\,Q^{\pi}\left( s,a \right),
\end{align}
If a policy $\pi^*$ satisfies $\mathrm{supp(}\pi ^*(\cdot |s))\subseteq \,\underset{a \in \calA}{\text{argmax}}\,\,Q^*\left( s,a \right)$ for any state $s\in \calS$, then it maximizes the state value function at all states simultaneously,
$$
V^{\pi^*}\left( s \right) =V^*\left( s \right),\quad \forall s\in \calS.
$$
Such policies are also referred to as the \textit{optimal policies}. Here, we use $\text{supp}(\cdot)$ to denote the support set of any probability distribution. 

\paragraph{Other useful notations in the proofs.}
Given a state $s \in \calS$, we use $\calA_s^\ast$ to denote the set of optimal actions at state $s$,
\begin{align*}
    \calA_s^\ast = \arg\max_{a\in\calA} Q^\ast(s,a).
\end{align*}
 Given a policy $\pi$,
we use $b_s^\pi$ to denote the probability of $\pi$ on non-optimal actions at state $s$,
\begin{align*}
    b_s^{\pi} = \sum_{a\notin \calA_s^{\ast}} \pi(a|s).
\end{align*}
When $\calA^*_s \ne \calA$, we also define the following two quantities:
\begin{align}
\widehat{V}^*\left( s \right) :=\underset{a\notin \calA _{s}^{*}}{\min}\,Q^*\left( s,a \right) ,\quad \widetilde{V}^*\left( s \right) :=\underset{a\notin \calA _{s}^{*}}{\max}\,Q^*\left( s,a \right),
\end{align}
We also define the optimal advantage function gap $\Delta$ as 
\begin{align}
\Delta :=\underset{s\in \tilde{S},a \notin \mathcal{A} _{s}^{*}}{\min}\left| A^{\pi^*}(s,a) \right|=\underset{s\in \tilde{S}}{\min} \, \{V^*(s)-\tilde{V}^*(s)\},
\label{optimal gap}
\end{align}
where $\widetilde{S}=\{s \in \calS:\calA^*_s \ne \calA\}$ denotes the set of states that have non-optimal actions.
The well-known Bellman optimality equation states that
\begin{align*}
V^*\left( s \right) =\underset{a\in \mathcal{A}}{\max}\,Q^*\left( s,a \right).
\end{align*}
Thus for any $\rho \in \Delta(\calA)$ that supported on $\calA \setminus \calA^*_s$ , it can be verified that
\begin{align*}
\widehat{V}^
\ast\left( s \right) \le \mathbb{E} _{a'\sim \rho}\left\{ Q^*\left( s,a' \right) \right\} \le \widetilde{V}^\ast\left( s \right).
\end{align*}
When $\calA ^\ast_s = \calA$, we let $\widehat{V}^\ast(s)= 0$ and $\widetilde{V}^\ast(s) = V^*(s)$. \\

The rest of this paper is organized as follows. In Section \ref{algorithm}, we derive the expression of the PG algorithm under the Hadamard parametrization, and the main theoretical results are provided in Section \ref{results}. In Section \ref{equivalent_algorithm}, we introduce  an equivalent expression of the algorithm based on the normalized Hadamard parametrization which does not require explicit projection onto the unit sphere. The proofs of the convergence results are provided in Section \ref{proofs} while the proofs of the useful lemmas are deferred to Section \ref{proofLemma}. Finally, we conclude this paper with a short discussion in Section \ref{conclusion}.
\section{Main results}

\subsection{Algorithm}
\label{algorithm}
Recall that  in the definition of the Hadamard parameterization (see \eqref{Hadamard Policy}), the vector $\theta_s$ for every state $s \in \calS$ is restricted to the unit sphere $\mathbb{S}^{|\calA|-1}$.  Thus the Riemannian gradient ascent algorithm can be used to solve the corresponding policy optimization problem.   To compute the Riemannian gradient, we can first compute the Euclidean gradient  and then projection it onto the tangent space. Let $g^E =[
    (g_{1}^E)^\tran,\cdots, (g_{|\calS|}^E)^\tran
]^\tran\in\R^{|\calS|\cdot |\calA|}$, 
where $g_s^E:=\nabla _{\theta _s}V^{\pi _{\theta}}\left( \mu \right) \in \mathbb{R} ^{|\calA |}$ is the Euclidean gradient with respect to parameter $\theta_s$, and let $g_{s,a}^E$ be the $a$-th entry of $g_s^E$. According to the policy gradient theorem \cite{Sutton1998}, the Euclidean gradient of the objective (\ref{PG Objective}) associated with any policy $\pi$ is given by
\begin{align*}
\frac{\partial V^{\pi}\left( \mu \right)}{\partial \pi \left( a|s \right)}=\frac{1}{1-\gamma}d_{\mu}^{\pi}\left( s \right) Q^{\pi}\left( s,a \right).
\end{align*}
Then it follows from the chain rule that
\begin{align*}
g_{s,a}^E=\frac{\partial V^{\pi _{\theta}}\left( \mu \right)}{\partial \theta _{s,a}}=\frac{\partial V^{\pi _{\theta}}\left( \mu \right)}{\partial \pi _{\theta}\left( a|s \right)}\cdot \frac{\partial \pi _{\theta}\left( a|s \right)}{\partial \theta _{s,a}}=\frac{2\theta _{s,a}}{1-\gamma}d_{\mu}^{\pi _{\theta}}\left( s \right) Q^{\pi}\left( s,a \right).
\end{align*}
Noting that the tangent space of $\mathbb{S}^{|\calA|-1}$  at each $\theta_s$ is given by 
\begin{align*}
T_{\theta _s}\mathbb{S} ^{|\calA |-1}=\left\{ y\in \mathbb{R} ^{|\calA |}: \left< y,\theta _s \right> =0 \right\},
\end{align*}
the projection of $g^E_s$ onto $T_{\theta _s}\mathbb{S} ^{|\calA |-1}$ (i.e.  the Riemannian gradient, denoted $g_s$) is given by 
\begin{align*}
g_s = (I-\theta_s\theta_s^\tran)g^E_s.
\end{align*}
Elementwisely, we have 
\begin{align*}
    g_{s,a} &= g_{s,a}^E-\theta_{s,a}(\theta_s^\tran g_s^E)\\
    &=g_{s,a}^E-\theta_{s,a}\sum_{a'}\theta_{s,a'}g^E_{s,a'}\\
    &=\frac{2\theta _{s,a}}{1-\gamma}d_{\mu}^{\pi _{\theta}}\left( s \right) Q^{\pi}\left( s,a \right)-\theta_{s,a}\sum_{a'}\frac{2\theta _{s,a'}^2}{1-\gamma}d_{\mu}^{\pi _{\theta}}\left( s \right) Q^{\pi}\left( s,a' \right)\\
    &=\frac{2\theta _{s,a}}{1-\gamma}d_{\mu}^{\pi _{\theta}}\left( s \right)\left(Q^\pi(s,a)-\sum_{a'}\theta_{s,a'}^2Q^\pi(s,a')\right)\\
    &=\frac{2\theta _{s,a}}{1-\gamma}d_{\mu}^{\pi _{\theta}}(s)A^\pi(s,a).\numberthis\label{Remannian Gradient}
\end{align*}



After updating $\theta_s$ towards the Riemannian gradient ascent direction, an additional normalization step is needed to ensure that the new estimate will be on the unit sphere. More precisely, letting $\eta$ be the step size, the basic update can be written as
\begin{align}
\theta _s\gets \frac{\theta _s+\eta g_s}{\left\| \theta _s+\eta g_s \right\| _2}.
\label{original update rule}
\end{align}
Substituting the expression of $g_{s,a}$ into  \eqref{original update rule} yields that
\begin{align}
\theta _{s,a}\gets \frac{\theta _{s,a}\left( 1+\frac{2\eta d_{\mu}^{\pi}\left( s \right)}{1-\gamma}A^{\pi}\left( s,a \right) \right)}{\sqrt{1+\frac{4\eta ^2\left( d_{\mu}^{\pi}\left( s \right) \right) ^2}{\left( 1-\gamma \right) ^2}\mathbb{E} _{a'\sim \pi \left( \cdot |s \right)}\left\{ \left( A^{\pi}\left( s,a' \right) \right) ^2 \right\}}}.
\label{update rule}
\end{align}
The complete description of the algorithm is presented in Algorithm~\ref{Hadamard PG}, where we simplify the notations $d^{\pi_{\theta^{k}}}_\mu$, $V^{\pi_{\theta^{k}}}$, $Q^{\pi_{\theta^{k}}}$, $A^{\pi_{\theta^{k}}}$,$\pi_{\theta^{k}}$ to $d^{k}_\mu$, $V^{k}$, $Q^{k}$, $A^{k}$, $\pi^{k}$, respectively.
Note that all the  entries of the initial parameters $\theta^0$ must be non-zeros, otherwise the parameters will not be updated. 

\begin{algorithm}[ht!]
\caption{Exact Policy Gradient Under Hadamard Parametrization}
\KwIn{step size $\eta > 0$, initial parameters $\{\theta^{0}_{s,a} \ne 0: s \in \calS, a\in \calA\}$}
\For{$k=0,1,2,...$}
{
    Compute the policy probabilities:
    $
 \pi _{s,a}^{k}=\left( \theta _{s,a}^{k} \right) ^2, ~\forall s\in \calS ,a\in \calA.
    $
    
    Compute $d^{k}_\mu$ and $A^{k}$ and update the parameter: 
    \begin{align}
\theta _{s,a}^{k+1}=\frac{\theta _{s,a}^{k}\left( 1+\frac{2\eta d_{\mu}^{k}\left( s \right)}{1-\gamma}A^k\left( s,a \right) \right)}{\sqrt{1+\frac{4\eta ^2\left( d_{\mu}^{k}\left( s \right) \right) ^2}{\left( 1-\gamma \right) ^2}\mathbb{E} _{a'\sim \pi ^k\left( \cdot |s \right)}\left\{ \left( A^k\left( s,a' \right) \right) ^2 \right\}}}, \quad \forall s\in \calS ,a\in \calA;
    \end{align}
   
    or equivalently, update the policy:
    \begin{align}
    \pi ^{k+1}(a|s)=\frac{\pi ^k\left( a|s \right) \left( 1+\frac{2\eta d_{\mu}^{k}\left( s \right)}{1-\gamma}A^k\left( s,a \right) \right) ^2}{1+\frac{4\eta ^2\left( d_{\mu}^{k}\left( s \right) \right) ^2}{\left( 1-\gamma \right) ^2}\mathbb{E} _{a'\sim \pi ^k\left( \cdot |s \right)}\left\{ \left( A^k\left( s,a' \right) \right) ^2 \right\}}.
    \end{align}

}
\label{Hadamard PG}
\end{algorithm}

\subsection{Theoretical results}
\label{results}

In this section, we provide the convergence results for Algorithm \ref{Hadamard PG} 
under standard assumptions. In Theorem \ref{Theorem1}, the global convergence of the algorithm is established. Theorem \ref{Theorem2} shows that the state value error decreases at an $O(\frac{1}{k})$ rate for all the iterations. Furthermore, we show that in Theorem \ref{Theorem3} that the algorithm enjoys a linear convergence rate after $k_0$ iterations and overall  displays a global linear convergence.

\begin{assumption}[Bounded Reward]
$r(s,a)\in[0,1],\forall (s,a).$
\label{assumption1}
\end{assumption}
\begin{assumption}[Sufficient Exploration]
$\tilde{\mu}:=\underset{s\in \calS}{\min}\,\mu \left( s \right) >0$.
\label{assumption2}
\end{assumption}
\begin{remark}
The existing global convergence results of the PG methods relies heavily on this assumption. 
Whether gradient ascent will globally converge if this condition is not met remains an open problem  \cite{pmlr-v125-agarwal20a}. When this condition is met, it follows from \eqref{d_eq} that the state visitation probability $d^\pi_\mu(s)$ is strictly positive for any $s \in \calS$.
\end{remark}
\begin{assumption} 
There exists at least one state that has non-optimal actions, i.e.  $\widetilde{S} \ne \emptyset$. Otherwise, it is trivial that any policy is an optimal policy.
\label{assumption3}


\end{assumption}

\begin{theorem}[Global Convergence]
\label{main result global}
For any constant step size  $\eta = \frac{(1-\gamma)^2}{4}\kappa$ with $\kappa \in (0,2)$ and the initial parameters $\theta^0$ such that $
\theta _{s,a}^{0} \ne 0 $ for all $(s,a)$ pairs, the value functions of the  policies generated by Algorithm~\ref{Hadamard PG}  converge  globally to $V^\ast(\rho)$ for any $\rho \in \Delta(\calS)$, i.e. 
\label{Theorem1}
\begin{align}
\underset{k\rightarrow \infty}{\lim}\,\,V^{k}\left( \rho \right) =V^*\left( \rho \right).
\end{align}
\end{theorem}

\begin{theorem}[Sublinear Convergence]
\label{main result sublinear}
Under the  conditions of Theorem \ref{main result global}, the value functions of the  policies generated by Algorithm \ref{Hadamard PG} converge to the optimal value with a sublinear rate,  
\label{Theorem2}
\begin{align}
V^*\left( \mu \right) -V^{k}\left( \mu \right) \le \frac{g\left( \kappa ,\lambda \right) ^{-1}}{k},
\label{SublinearRate}
\end{align}
where $
g\left( \kappa ,\lambda \right) =\frac{2\kappa \tilde{\mu}^2\left( 1-\gamma \right)^4 \lambda}{4+\kappa ^2}$ and $\lambda=\inf_{k\geq 0} \left\{ \min_{s\in\calS} (1-b_s^k) \right\}$ (here $b_s^k$ is short for $b_s^{\pi_k}$) is a strictly positive constant (see Lemma \ref{lemLambda}).
\end{theorem}

\begin{theorem}[Linear Convergence]
\label{main result linear}
For any $k\ge 0$ and $\delta \in (0,1)$, under the  conditions of Theorem \ref{main result global}, the value functions of the  policies generated by Algorithm \ref{Hadamard PG} satisfy
\begin{align*}
V^*\left( \mu \right) -V^k\left( \mu \right) \le \frac{1}{(1-\gamma)^2}\cdot \left( 1-\rho \left(\delta, \kappa ,\gamma \right) \right) ^{k-k_0\left( \delta, \kappa ,\gamma \right)},
\end{align*}
where $\rho \left(\delta, \kappa ,\gamma \right) = \frac{3\kappa^2 (1-\gamma)^2 \tilde{\mu}}{16}\cdot \Delta\cdot(1-\delta)$ and $k_0\left( \delta, \kappa ,\gamma \right) =  \frac{g(\kappa,\lambda)^{-1}}{\tilde{\mu} \Delta}\cdot\min\{\frac{3\delta\kappa}{20\gamma},\frac{1}{2}\}^{-1}$
\label{Theorem3}
\end{theorem}

\begin{remark}\kw{
It can be observed that there are two problem-dependent parameters $\tilde{\mu}$ and $\Delta$ in $\rho \left(\delta, \kappa ,\gamma \right)$. Firstly, the parameter $\tilde{\mu}$ essentially comes from the lower bound for the visitation measure term $d_\mu^{\pi_\theta}(s)$ in the Riemannian gradient expression \eqref{Remannian Gradient}. Due to the existence of this term, the effective step size is indeed $\eta \cdot\min_s d_\mu^{\pi_\theta}(s)\propto \eta\cdot\tilde{\mu}$ (since $(1-\gamma)\tilde{\mu}$ provides a uniform lower bound for $d_\mu^{\pi_\theta}(s))$, which means smaller $\tilde{\mu}$ leads to smaller effective step size, and thus slower convergence. Secondly, recall that $\Delta$ is the optimal advantage function gap, and small $\Delta$ means that it is difficult to tell the optimal actions apart from non-optimal ones. Thus, it is reasonable for an algorithm to have a slow convergence rate when $
\Delta$ is small.
}\end{remark}

\subsection{Equivalent algorithm without spherical constraint}
\label{equivalent_algorithm}

\begin{algorithm}[ht!]
\caption{Deterministic Policy Gradient Under Normalized Hadamard Parametrization}
\label{Normalized Hadamard PG}
\KwIn{step size $\eta > 0$, initial parameters $\{\theta^{0}_{s,a} \ne 0 : s \in \calS, a\in \calA\}$}
\For{$k=0,1,2,...$}
{
    Compute the policy probabilities:
    \begin{align*}
\forall s\in \calS ,a\in \calA :\,\pi _{s,a}^{k}=\frac{\left( \theta _{s,a}^{k} \right) ^2}{\left\| \theta _{s}^{k} \right\| _{2}^{2}}.
    \end{align*}
    \\
    Compute the exact gradient of $\mathcal{L} _k\left( \theta \right)$ in \eqref{eq:tmp_ke01} using \eqref{normalized hadamard gradient}.
    \\
    Update the parameter:
\begin{align*}
\theta _{k+1}=\theta _k+ \eta \cdot \nabla _{\theta}\mathcal{L} _k\left( \theta \right) |_{\theta =\theta _k}.
\end{align*}
}
\end{algorithm}

Note that it requires explicit normalization to satisfy the spherical constraint in Algorithm~\ref{Hadamard PG}, which is not flexible when using function approximation. In this section, we are about to show that Algorithm~\ref{Hadamard PG} is equivalent to an algorithm that does not require parameter normalization. To this end, first define the {normalized Hadamard parametrization}  as
\begin{align}
\pi _{\theta}\left( a|s \right) =\frac{(\theta _{s,a})^2}{\left\| \theta _s \right\|^2_2}, \quad\theta _s\in \mathbb{R} ^{|\calA |}\setminus \left\{ \mathbf{0} \right\}.
\label{normalized hadamard}
\end{align}
which  can be seen as a special case of escort policy with $p =2$ \cite{Mei_Xiao_Dai_Li_Szepesvari_Schuurmans_2020}. In the $k$-th iteration, the algorithm optimizes following objective function:
\begin{align}\label{eq:tmp_ke01}
\mathcal{L} _k\left( \theta \right) :=\frac{1}{1-\gamma}\mathbb{E} _{s\sim d_{\mu}^{k}}\{ \left\| \theta _{s}^{k} \right\| _{2}^{2} \cdot \mathbb{E} _{a\sim \pi _{\theta}}\{ A^k\left( s,a \right) \} \},
\end{align}
and a gradient ascent step is conducted  with a constant step size $\eta$,
\begin{align}
\theta _{k+1}=\theta _k+ \eta \cdot \nabla _{\theta}\mathcal{L} _k\left( \theta \right) |_{\theta =\theta _k}.
\label{normalized hadamard update rule}
\end{align}
For any $s\in \calS$, a direct computation yields that
\begin{align}
\left. \frac{\partial \mathcal{L} _k\left( \theta \right)}{\partial \theta _{s,a}} \right|_{\theta =\theta ^k}&=\sum_{a'\in \mathcal{A}}{\left. \frac{\partial \mathcal{L} _k\left( \theta \right)}{\partial \pi \left( a'|s \right)} \right|_{\theta =\theta ^k}\cdot \left. \frac{\partial \pi \left( a'|s \right)}{\partial \theta _{s,a}} \right|}_{\theta =\theta ^k} \notag
\\
\quad         &=\sum_{a'\in \mathcal{A}}{\frac{\left\| \theta _{s}^{k} \right\| _{2}^{2}d_{\mu}^{k}\left( s \right) A^k\left( s,a' \right)}{1-\gamma}\cdot \frac{2\theta _{s,a}^{k}}{\left\| \theta _{s}^{k} \right\| _{2}^{4}}\left( \mathbf{1}\left[ a=a' \right] \left\| \theta _{s}^{k} \right\| _{2}^{2}-\left( \theta _{s,a'}^{k} \right) ^2 \right)} \notag
\\
\quad         &=\frac{2d_{\mu}^{k}\left( s \right)}{1-\gamma}\frac{\theta _{s,a}^{k}}{\left\| \theta _s \right\| _{2}^{2}}\left( A^k\left( s,a \right) \left\| \theta _{s}^{k} \right\| _{2}^{2}-\mathbb{E} _{a'\sim \pi ^k\left( \cdot |s \right)}\left\{ A^k\left( s,a' \right) \right\} \right) \notag
\\
\quad         &=\frac{2d_{\mu}^{k}\left( s \right)}{1-\gamma}\cdot \theta _{s,a}^{k}\cdot A^k\left( s,a \right).
\label{normalized hadamard gradient}
\end{align}
The complete procedure is summarized in Algorithm \ref{Normalized Hadamard PG}.
Since the normalized Hadamard parametrization is not defined when $\theta_s = \textbf{0}$ , one must ensure that $ \left\| \theta _{s}^{k} \right\| _{2}^{2} \ne 0$ for all iterations and $\forall s \in \calS$. However, this is not an issue if $\left\| \theta _{s}^{0} \right\| _{2}^{2} \ne 0$ since
\begin{align*}
\left\|\theta _{s}^{k+1} \right\| _2=\left\| 
	\theta _{s}^{k}+\eta \cdot \left. \nabla _{\theta _s}\mathcal{L} _k\left( \theta \right) \right|_{\theta =\theta _k} \right\| _2\overset{\left( a \right)}{=}\sqrt{\left\| \theta _{s}^{k} \right\| _{2}^{2}+\eta ^2\left\| \left. \nabla _{\theta _s}\mathcal{L} _k\left( \theta \right) \right|_{\theta =\theta _k} \right\| _{2}^{2}}\ge \left\| \theta _{s}^{k} \right\| _2,
\end{align*}
where step (a) is due to (\ref{normalized hadamard gradient is orthogonal}).

Next we show that Algorithm~\ref{Hadamard PG} and Algorithm~\ref{Normalized Hadamard PG} are equivalent to each other in the policy domain provided the initial policies are the same. First, 
it is easy to see that
\begin{align}
\left< \left. \nabla _{\theta _s}\mathcal{L} _k\left( \theta \right) \right|_{\theta =\theta ^k},\theta _{s}^{k} \right> =\sum_{a\in \mathcal{A}}{\frac{2d_{\mu}^{k}\left( s \right)}{1-\gamma}\left( \theta _{s,a}^{k} \right) ^2A^k\left( s,a \right)}=\frac{2d_{\mu}^{\pi}\left( s \right)}{1-\gamma}\left\| \theta _{s}^{k} \right\| _{2}^{2}\mathbb{E} _{a\sim \pi ^k\left( \cdot |s \right)}\left\{ A^k\left( s,a \right) \right\} =0.
\label{normalized hadamard gradient is orthogonal}
\end{align}
Plugging (\ref{normalized hadamard}), (\ref{normalized hadamard gradient}) and (\ref{normalized hadamard gradient is orthogonal}) into (\ref{normalized hadamard update rule}) yields that 
\begin{align*}
\pi ^{k+1}(a|s)&=\left( \frac{\theta _{s,a}^{k+1}}{\left\| \theta _{s}^{k+1} \right\|_2} \right) ^2
\\
&=\frac{\left( \theta _{s,a}^{k}+\left. \eta \cdot \frac{\partial \mathcal{L} _k\left( \theta \right)}{\partial \theta _{s,a}} \right|_{\theta =\theta ^k} \right) ^2}{\left\| \theta _{s}^{k}+\left. \eta \cdot \nabla _{\theta _s}\mathcal{L} _k\left( \theta \right) \right|_{\theta =\theta ^k} \right\| _{2}^{2}}
\\
&=\frac{\left( \theta _{s,a}^{k}\left( 1+\frac{2\eta d_{\mu}^{k}\left( s \right)}{1-\gamma}A^k\left( s,a \right) \right) \right) ^2}{\sum_{a'}{\left( \theta _{s,a'}^{k} \right) ^2\left( 1+\frac{2\eta d_{\mu}^{k}\left( s \right)}{1-\gamma}A^k\left( s,a' \right) \right) ^2}}
\\
&=\frac{\frac{\left( \theta _{s,a}^{k} \right) ^2}{\left\| \theta _{s}^{k} \right\| _{2}^{2}}\cdot \left( 1+\frac{2\eta d_{\mu}^{k}\left( s \right)}{1-\gamma}A^k\left( s,a \right) \right) ^2}{\sum_{a'}{\frac{\left( \theta _{s,a'}^{k} \right) ^2}{\left\| \theta _{s}^{k} \right\| _{2}^{2}}\left( 1+\frac{2\eta d_{\mu}^{k}\left( s \right)}{1-\gamma}A^k\left( s,a' \right) \right) ^2}}
\\
&=\frac{\pi ^k\left( a|s \right) \left( 1+\frac{2\eta d_{\mu}^{k}\left( s \right)}{1-\gamma}A^k\left( s,a \right) \right) ^2}{1+\frac{4\eta ^2\left( d_{\mu}^{k}\left( s \right) \right) ^2}{\left( 1-\gamma \right) ^2}\mathbb{E} _{a'\sim \pi ^k\left( \cdot |s \right)}\left\{ \left( A^k\left( s,a' \right) \right) ^2 \right\}},
\end{align*}
which is the same as the policy update in Algorithm \ref{Hadamard PG}. Therefore, both algorithms output the same polices at all the iterations if using the same initial policy $\pi^{0}$. 

\kw{Additionally, it can be directly verified that the policy gradient direction under the normalized Hadadmard parameterization (also known as escort policy wiht $p=2$ in \cite{Mei_Xiao_Dai_Li_Szepesvari_Schuurmans_2020}) 
is given by 
\begin{align*}
\frac{2d_{\mu}^{k}\left( s \right)}{1-\gamma}\cdot \theta _{s,a}^{k}\cdot A^k\left( s,a \right)\frac{1}{\|\theta_s^k\|_2^2},
\end{align*}
which only differs from \eqref{normalized hadamard gradient} by a state dependent factor $\frac{1}{\|\theta_s^k\|_2^2}$.
Thus Algorithm \ref{Normalized Hadamard PG} can be seen as escort PG with preconditioning or with state-dependent step sizes $\eta_k(s) = \eta \cdot \left\|\theta_s^k \right\|^2_2$.}

\subsection{Numerical simulations}
\label{simulation}
In this section, we conduct numerical experiments to (1) demonstrate the influences of the parameters on the convergence of Hadamard PG, and (2) compare the performance of Hadamard PG with  other typical algorithms: PG under direct parameterization and softmax parameterization, and NPG under softmax parameterization. \kw{The numerical experiments are conducted on random MDPs, in which the rewards for the state-action pairs and the entries of the transition probability matrix are generated through a uniform distribution over $\left(0,1\right)$.} We then re-scale each row of the transition probability matrix to obtain a stochastic matrix. In the simulations, we use $\gamma = 0.9$ and $\mu(s)=\frac{1}{\left|\calS\right|}$.

\begin{figure}[ht!]
    \centering
    \subfigure[]{
    \includegraphics[width=0.48\hsize, height=0.36\hsize]{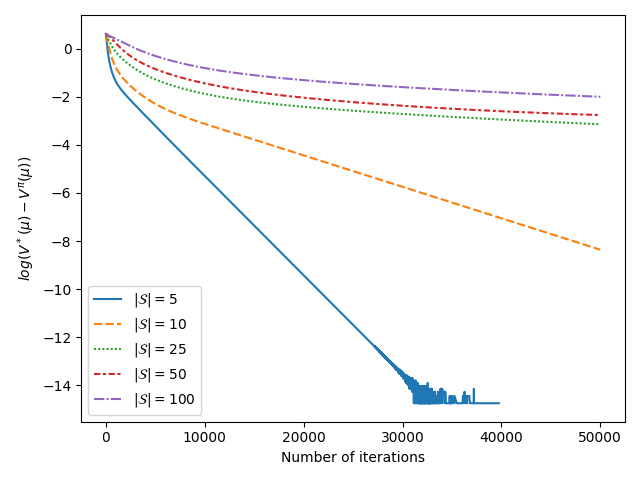}
    \label{fig: sub_figure1}
    }
    \subfigure[]{\includegraphics[width=0.48\hsize, height=0.36\hsize]{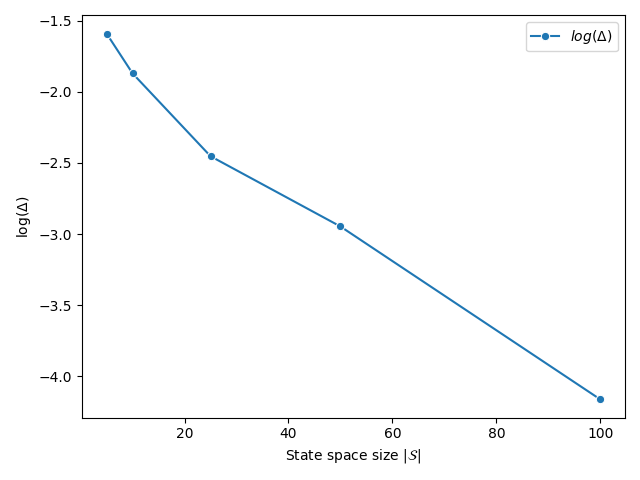}\label{fig: sub_figure2}}
    \caption{{(a) Convergence plot of Hadamard PG on random MDPs; (b) Plot of $\Delta$ versus $|\calS|$.}}
    \label{fig: figure1}
\end{figure}

\begin{figure}[ht!]
    \centering
    \subfigure[Bandit]{
    \includegraphics[width=0.48\hsize, height=0.36\hsize]{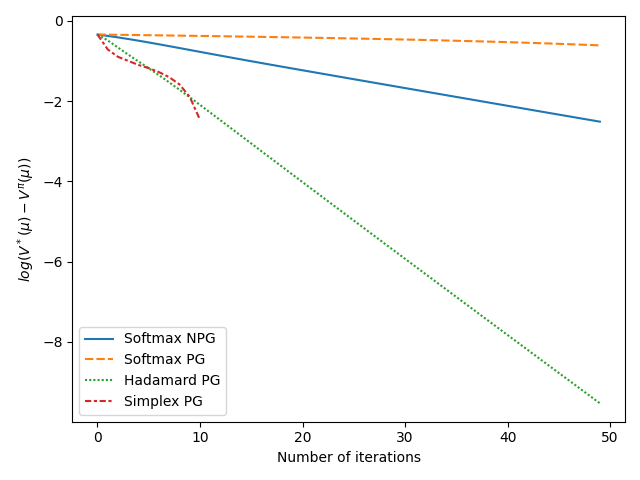}
    \label{fig2: sub_figure_bandit}
    }
    \subfigure[MDP]{\includegraphics[width=0.48\hsize, height=0.36\hsize]{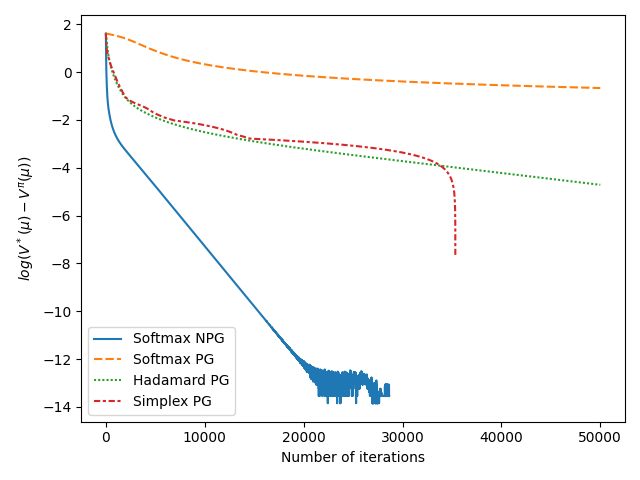}\label{fig2: sub_figure_MDP}}
    \caption{{Performance comparison between algorithms on a random bandit problem with $\left|\calA\right|$ = 10 and a  random MDP with $\left|\calS\right|$ = 50, $\left|\calA\right|$ = 10}.}
    \label{fig2: figure1}
\end{figure}

\paragraph{Convergence of Hadamard PG with different $\left|\calS\right|$.} Here we \kw{empirically} demonstrate the influences of $\tilde{\mu}$ and $\Delta$ in the theoretical linear convergence rate presented in Theorem \ref{main result linear}. The step size in Hadamard PG is chosen to be $\eta = \frac{(1-\gamma)^2}{2}$. Since $\mu(s)=\frac{1}{|\mathcal{S}|}$, we have $\tilde{\mu}=\frac{1}{|\mathcal{S}|}.$ For the term $\Delta$, larger problem size (i.e. the larger $\left|\calS\right|$ and $\left|\calA\right|$) often induces  smaller $\Delta$. In our simulations, we fix $|\calA|=10$ and vary $|\calS|$ from 5 to 100. The convergence plots of Hadamard PG under different $|\calS|$ are displayed in Figure~\ref{fig: sub_figure1}. It is evident that the convergence of Hadamard PG becomes slower as $|\calS|$ increases. In addition, the values of $\Delta$ versus $|\calS|$ are plotted in Figure~\ref{fig: figure1}, which confirms that $\Delta$ decreases as $|\calS|$ increases.
\paragraph{Comparison with other algorithms.} In this part, we compare the performance of Hadamard PG with the other \kw{aforementioned} algorithms on a random bandit problem and a random MDP problem. Bandit is a special RL problem  with one single state and $\gamma = 0$.   The update rules of the other test algorithms are as follows:
\begin{itemize}
    \item Simplex PG:
\begin{align*}
\pi ^{k+1}\left( \cdot |s \right) =\text{Proj}_{\Delta \left( \mathcal{A} \right)}\left( \pi ^k\left( \cdot |s \right) +\frac{\eta d_{\mu}^{k}\left( s \right)}{1-\gamma}Q^k\left( s,a \right) \right)
\end{align*}
where  $\text{Proj}_{\Delta(\mathcal{A})}$ denotes the projection onto $\Delta(\mathcal{A})$ in Euclidean norm, i.e.  $\mathrm{Proj}_{\Delta \left( \mathcal{A} \right)}\left( v \right) =\underset{p\in \Delta \left( \mathcal{A} \right)}{\mathrm{arg}\min}\,\,\left\| p-v \right\| _{2}^{2}$.
    
    \item Softmax PG:
 \begin{align*}
\theta _{s,a}^{k+1}=\theta _{s,a}^{k}+\frac{\eta d_{\mu}^{k}\left( s \right)}{1-\gamma}\cdot \pi ^k\left( a|s \right) \cdot A^k\left( s,a \right). 
 \end{align*}
    \item Softmax NPG:
\begin{align*}
\theta _{s,a}^{k+1}=\theta _{s,a}^{k}+\frac{\eta}{1-\gamma}\cdot A^k\left( s,a \right).
 \end{align*}
\end{itemize}
For fair comparison, we use $\eta = \frac{1-\gamma}{2}$ for all algorithms (\kw{it is larger than the upper bound $\frac{(1-\gamma)^2}{2}$ used in our theoretical analysis which suggests the possibility to extend the analysis to a wider range of step sizes}). Figure \ref{fig2: sub_figure_bandit} shows the performance comparison on the random bandit problem with $\left|\calA\right|$ = 10. It can be observed that Hadamard PG has a better performance  for this problem. Figure \ref{fig2: sub_figure_MDP} shows the convergence of different algorithms on the random MDP with $\left|\calS\right|$ = 50 and $\left|\calA\right|$ = 10, \kw{in which case softmax NPG performs best}. A \kw{potential} explanation is that the step size in softmax NPG is \textbf{not} slowed down by the state visitation measure $d^k_\mu(s)$ which is often significantly small when the state space size $\left|\calS\right|$ is large. 
This also suggests that combining the idea of NPG with Hadamard parameterization
may yield an algorithm with better convergence, which will be left as future work.

\section{Proofs}
\label{proofs}
We first list several useful lemmas. The first one is the  performance difference lemma whose proof can be found in  \cite{kakade2002approximately} and the proofs of the other lemmas are deferred to Section \ref{proofLemma}. 

\begin{lemma}[Performance difference lemma \cite{kakade2002approximately}]
\label{performanceDifferenceLemma}
	For any two policy $\pi_1, \pi_2$, and any $\rho \in \Delta(\calS)$, one has
	\begin{align*}
		V^{\pi_1}(\rho) - V^{\pi_2}(\rho) =\frac{1}{1-\gamma} \E[s\sim d_{\rho}^{\pi_1}]{ \E[a \sim \pi_1(\cdot|s)]{A^{\pi_2}(s,a)} }.
	\end{align*}
\end{lemma}


\begin{lemma}
\label{keyLemma1}
For any policy $\pi$ and $\rho \in \Delta(\calS)$,
\begin{align}
\Delta \cdot \mathbb{E} _{s\sim \rho}\{b^\pi_s\} \le 
V^\ast(\rho)-V^\pi(\rho) \le 
\frac{1}{(1-\gamma)^2}\mathbb{E} _{s\sim d^\pi_\rho}\{b^\pi_s\}.
\end{align}
\end{lemma}

\begin{lemma}[One-iteration improvement]
\label{improvedLowerBound}
Consider the policy $\pi^k$ generated by Algorithm \ref{Hadamard PG}. Suppose $\eta = \frac{(1-\gamma)^2\kappa}{4}$ for $\kappa \in (0,1)$. Then for any $\rho \in \Delta(\calS)$, one has
\begin{align*}
V^{k+1}(\rho )-V^k(\rho )\ge \frac{2\kappa \tilde{\mu}\left( 1-\gamma \right) ^2}{4+\kappa ^2} \cdot \underset{s\in \mathrm{supp}\left( \rho \right)}{\min}\left\{ \rho \left( s \right) \right\} \cdot \sum_{s\in \mathrm{supp}\left( \rho \right)}{\sum_{a \in \calA}{\pi ^k}(a|s)(A^k(s,a))^2}\ge 0.
\end{align*}
In particular when $\rho = \mu$,
\begin{align}
V^{k+1}(\mu )-V^k(\mu )\ge \frac{2\kappa \tilde{\mu}^2\left( 1-\gamma \right) ^2}{4+\kappa ^2}\sum_s{\sum_a{\pi ^k}(a|s)(A^k(s,a))^2}\ge 0.
\end{align}
\end{lemma}
\begin{lemma}
\label{policyUpdate}
 Consider the policy $\pi^k$ generated by Algorithm \ref{Hadamard PG}. One has
	\begin{align}
		\pi^{k+1}(a|s) - \pi^k(a|s) =\frac{ \pi^k(a|s) }{1+ \eta^2 \twonorm{\vg_s^k}^2} \cdot \frac{4\eta d_{\mu}^k(s)}{1-\gamma} \left( A^k(s,a) + \frac{\eta d_{\mu}^k(s)}{1-\gamma} \left(  (A^{k}(s,a))^2 - \E[a'\sim \pi^k(\cdot|s)]{(A^{k}(s,a'))^2}\right)\right). 
	\end{align}
\end{lemma}

\begin{lemma}
	\label{lemma b}
        For any $s \in \calS$, the sequence $\{b_s^k\}_{k\geq k_0}$ generated by Algorithm \ref{Hadamard PG} satisfies 
        \begin{align*}
        b_{s}^{k}-b_{s}^{k+1}\ge \frac{\eta d_{\mu}^{k}(s)}{2(1-\gamma)}\cdot \frac{\min\mathrm{(}1-b_{s}^{k},b_{s}^{k})}{1+\eta ^2\left\| g_{s}^{k} \right\| _{2}^{2}}\cdot \left(3\kappa\Delta -20\gamma \left\|V^*-V^k\right\|_\infty \right),
        \end{align*}
\end{lemma}

\begin{lemma}
	\label{lemLambda}
	Suppose $\pi^0(a|s)>0$ for all $s\in\calS, a\in\calA$. Then $\forall s,a,k$, one has $\pi^k(a|s) >0$. Moreover,
	\begin{align}
		\label{eqLambda}
		\lambda:=\inf_{k\geq 0} \left\{ \min_{s\in\calS} (1-b_s^k) \right\} >0.
	\end{align}
\end{lemma}

\subsection{Proof of Theorem \ref{main result global}}
\label{proof them1}
\kw{The argument  for the global convergence of Hadamard PG is overall similar to that for  Theorem 5.1 in \cite{Agarwal_Kakade_Lee_Mahajan_2019}.}
We first show that the value function improves monotonically and then it is easy to further show that the limit of advantage function is non-positive, which directly implies the global convergence. 
\begin{proof}[Proof of Theorem \ref{main result global}]
Firstly,  Lemma \ref{improvedLowerBound} implies that $\{V^k(\rho)\}_k$ is nondecreasing. Since $V^k(\rho)$ is upper bounded, we conclude that $\{V^k(\rho)\}$ is convergent. Similarly, we can also prove that the sequence $\{Q^k(s,a)\}_k$ is convergent for any $(s,a)$. Specifically, one has
\begin{align*}
	Q^{k+1}(s,a) = r(s,a) + \E[s'\sim P(\cdot|s,a)]{V^{k+1}(s')}\geq r(s,a) +  \E[s'\sim P(\cdot|s,a)]{V^{k}(s')} = Q^k(s,a),
\end{align*}
where the inequality  is due to the fact that $V^k(\rho)$ is nondecreasing with $\rho=P(\cdot|s,a)$. Since $Q^k(s,a)$ is upper bounded for any $(s,a)$, it implies that the sequence is convergent. Define
\begin{align*}
	V^{\infty}(\rho) &= \lim_{k\rightarrow +\infty} V^k(\rho),\\
	Q^{\infty}(s,a)  &=\lim_{k\rightarrow +\infty} Q^{k}(s,a).
\end{align*}
Moreover, define
\begin{align*}
	A^{\infty}(s,a)=\lim_{k\rightarrow +\infty}A^k(s,a) = \lim_{k\rightarrow +\infty}Q^{k}(s,a) -\lim_{k\rightarrow +\infty} V^k(s).
\end{align*}
We  will show that 
\begin{align}
	\label{aux1}
	A^{\infty}(s,a) \leq 0
\end{align}
for any $(s,a)$. According to Lemma \ref{improvedLowerBound}, we have
\begin{align*}
	V^{k+1}(\mu) - V^{k}(\mu )\geq  \frac{3\kappa \tilde{\mu}^2  (1-\gamma)^2 }{ (4+ \kappa^2 )}   \sum_s  \sum_a \pi^{k}(a|s) (A^{k}(s,a))^2\geq 0.
\end{align*}
Taking a limit on both sides yields that
\begin{align*}
	0=\lim_{k\rightarrow +\infty} V^{k+1}(\mu) - V^{k}(\mu ) \geq \frac{3\kappa \tilde{\mu}^2  (1-\gamma)^2 }{ (4+ \kappa^2 )} \lim_{k\rightarrow +\infty}  \sum_s  \sum_a \pi^{k}(a|s) (A^{k}(s,a))^2\geq 0,
\end{align*}
which implies that 
\begin{align}
	\label{aux2}
	\lim_{k\rightarrow +\infty} \E[a\sim \pi^k(\cdot|s)]{(A^{k}(s,a))^2}=0
\end{align}
for any state $s\in\calS$. We will use proof by contradiction to show that $A^\infty(s,a)\leq 0$ for any state $s$ and any action $a$. Suppose that there exists a state-action pair $(s,a)$ such that $A^\infty(s,a)>0$. 
According to Lemma \ref{policyUpdate}, it can be seen that the sign of $\pi^{k+1}(a|s) - \pi^k(a|s)$ only depends on the term
\begin{align*}
	I_{s,a}^k:= A^k(s,a) + \frac{\eta d_{\mu}^k(s)}{1-\gamma} \left(  (A^{k}(s,a))^2 - \E[a'\sim \pi^k(\cdot|s)]{(A^{k}(s,a'))^2}\right)
\end{align*}
and we have
\begin{align*}
\liminf_{k\rightarrow +\infty}\,\,I_{s,a}^{k} &\ge \liminf_{k\rightarrow +\infty}\,\,A^k(s,a)+\eta \tilde{\mu}\,\liminf_{k\rightarrow +\infty}\left( (A^k(s,a))^2-\mathbb{E} _{a'\sim \pi ^k(\cdot |s)}(A^k(s,a'))^2 \right) 
\\
\quad            &=A^{\left( \infty \right)}\left( s,a \right) +\eta \tilde{\mu}\,(A^{\left( \infty \right)}(s,a))^2
\\
\quad            &>A^{\left( \infty \right)}\left( s,a \right), 
\end{align*}
where the second line follows from \eqref{aux1} and \eqref{aux2}. Thus it implies that there exists a constant $N_0$ such that 
\begin{align*}
	I_{s,a}^k > \frac{A^\infty(s,a)}{2} > 0
\end{align*}
for any $k\geq N_0$ and $\{\pi^k\}_{k\geq N_0}$ is a nondecreasing sequence. Since $\pi^{N_0}(a|s)>0$ (see the first claim of Lemma~\ref{lemLambda}), we have
\begin{align*}
	\lim_{k\rightarrow +\infty} \pi^k(a|s) >0 \text{ and } \lim_{k\rightarrow +\infty} \left( \pi^k(a|s)(A^k(s,a))^2 \right) >0,
\end{align*}
which contradicts with \eqref{aux2}. Therefore, $A^{\infty}(s,a)\leq 0$ for any $(s,a)$. Finally, we show that $V^{\infty}(\rho) = V^\ast(\rho)$. Using Lemma \ref{performanceDifferenceLemma}, we have
\begin{align*}
	V^\ast(\rho) - V^k(\rho) = \frac{1}{1-\gamma}\sum_{s\in \calS} d_{\rho}^\ast(s) \sum_{a\in\calA}\pi^\ast(a|s) A^k(s,a).
\end{align*}
Taking  a limit on both sides yields that
\begin{align*}
	V^\ast(\rho) - \lim_{k\rightarrow +\infty} V^k(\rho) = \frac{1}{1-\gamma}\sum_{s\in \calS} d_{\rho}^\ast(s) \sum_{a\in\calA}\lim_{k\rightarrow +\infty}\pi^\ast(a|s) A^k(s,a)\leq 0,
\end{align*}
where the inequality  is due to $A^{\infty}(s,a)\leq 0$. 
\end{proof}

\subsection{Proof of Theorem \ref{main result sublinear}}
\label{proofLemma46}
\kw{The proof of the sublinear convergence of Hadamard PG is substantially different from for example the proof for the sublinear convergence of simplex PG and softmax PG in \cite{Agarwal_Kakade_Lee_Mahajan_2019, linxiao2022, Mei_Xiao_Dai_Li_Szepesvari_Schuurmans_2020}. Instead of using the smooth property of the value function together with the gradient dominance property}, we directly establish a lower bound for the one step improvement $V^{k+1}-V^k$, \kw{as well as} an upper bound for the distance of $V^k$ to the optimal value  $V^*$. Together, we can obtain an inequality of the following form:
$$V^{k+1}-V^k \gtrsim (V^*-V^k)^2,$$ 
from which the sublinear convergence follows directly.
\begin{proof}[Proof of Theorem \ref{main result sublinear}]
First define the auxiliary policy supported on $\calA_s^*$ as follows:
\begin{align*}
\psi^k(\cdot|s) = \begin{cases}
\frac{\pi^k \left( \cdot |s \right)}{1-b_{s}^{k}}\mathbf{1}\left[ \cdot \in \mathcal{A} _{s}^{*} \right],\mbox{ if }b_s^k\neq 1,\\
\mbox{any distribution supported on }\calA_s^*,&\mbox{ if }b_s^k=1.
\end{cases}
\end{align*}
It is clear that $\pi^k(a|s)=(1-b_s^k)\psi^k(a|s)$ for $a\in\calA_s^*$. 
Moreover, 
since $\text{supp}(\psi^k) \subseteq \calA^*_s$, one can verify that $\psi^k$ is an optimal policy. Actually, by the performance difference lemma, 
\begin{align*}
	V^{\psi^k}(\mu) - V^\ast(\mu) =\frac{1}{1-\gamma}\E[s\sim d_{\mu}^{\psi^k}]{\sum_{a\in\calA_s^\ast} \psi^k(a|s) A^\ast(s,a) }=0,\numberthis\label{eq 5}
\end{align*}
where the last line is due to $A^\ast(s,a)=0$ for any $a\in\calA_s^\ast$.
A direct computation yields that
	\begin{align*}
		V^{k+1}(\mu) - V^{k}(\mu ) &\stackrel{(a)}{\geq} \frac{3\kappa \tilde{\mu}^2  (1-\gamma)^2 }{ (4+ \kappa^2 )}   \sum_s  \sum_a \pi^{k}(a|s) (A^{k}(s,a))^2\\
		&\geq \frac{2\kappa \tilde{\mu}^2  (1-\gamma)^2 }{ (4+ \kappa^2 )}   \sum_s d_{\mu}^{\psi^k}(s) \sum_a \pi^{k}(a|s) (A^{k}(s,a))^2\\
		&=\frac{2\kappa \tilde{\mu}^2  (1-\gamma)^2 }{ (4+ \kappa^2 )}   \sum_s d_{\mu}^{\psi^k}(s) \left( \sum_{a\in\calA_s^\ast} \pi^{k}(a|s) (A^{k}(s,a))^2 + \sum_{a\notin\calA_s^\ast} \pi^{k}(a|s) (A^{k}(s,a))^2 \right)\\
		&\geq \frac{2\kappa \tilde{\mu}^2  (1-\gamma)^2 }{ (4+ \kappa^2 )}   \sum_s d_{\mu}^{\psi^k}(s)   \sum_{a\in\calA_s^\ast} \pi^{k}(a|s) (A^{k}(s,a))^2 \\
		&=\frac{2\kappa \tilde{\mu}^2  (1-\gamma)^2 }{ (4+ \kappa^2 )}   \sum_s d_{\mu}^{\psi^k}(s) (1-b_s^k)  \sum_{a\in\calA } \psi^k(a|s) (A^{k}(s,a))^2 \\
		&\geq \frac{2\kappa \tilde{\mu}^2  (1-\gamma)^2 }{ (4+ \kappa^2 )}  \left( \min_s (1-b_s^k)\right)  \sum_s d_{\mu}^{\psi^k}(s)  \sum_{a\in\calA } \psi^k(a|s) (A^{k}(s,a))^2 \\
		&\stackrel{(b)}{\geq}\frac{2\kappa \tilde{\mu}^2  (1-\gamma)^2 \lambda }{ (4+ \kappa^2 )}  \sum_s d_{\mu}^{\psi^k}(s)  \sum_{a\in\calA } \psi^k(a|s) (A^{k}(s,a))^2\\
		&\stackrel{(c)}{\geq} \frac{2\kappa \tilde{\mu}^2  (1-\gamma)^2 \lambda }{ (4+ \kappa^2 )}  \left(\sum_s d_{\mu}^{\psi^k}(s)  \sum_{a\in\calA } \psi^k(a|s) |A^{k}(s,a)| \right)^2, \numberthis\label{eq46a}
	\end{align*}
	where step (a) is due to Lemma \ref{improvedLowerBound}, step (b) follows from Lemma \ref{lemLambda}, and step (c) is due to Jensen's inequality. On the other hand, we have
	\begin{align*}
		V^\ast(\mu) - V^k(\mu) &\stackrel{(a)}{=} V^{\psi^k}(\mu) - V^k(\mu)\\
		&\stackrel{(b)}{=} \frac{1}{1-\gamma}\E[s\sim d_{\mu}^{\psi^k}]{\E[a\sim \psi^k]{A^k(s,a)}}\\
		&\leq \frac{1}{1-\gamma}\E[s\sim d_{\mu}^{\psi^k}]{\E[a\sim \psi^k]{ \left| A^k(s,a) \right| }},\numberthis \label{eq46b}
	\end{align*}
	where step (a) is due to \eqref{eq 5}, step (b) is due to the performance difference lemma. Combining \eqref{eq46a} and \eqref{eq46b} together, we have
	\begin{align*}
		V^{k+1}(\mu) - V^{k}(\mu )  &\geq \frac{2\kappa \tilde{\mu}^2  (1-\gamma)^4\lambda }{ 4+ \kappa^2 }  \left(V^\ast(\mu) - V^k(\mu)  \right)^2\\
		:&=g(\kappa,\lambda)  \left(V^\ast(\mu) - V^k(\mu)  \right)^2.
	\end{align*}
	Let $\delta_k = V^\ast(\mu) - V^k(\mu)$. It can be seen that
	\begin{align*}
		\delta_k - \delta_{k+1} = V^{k+1}(\mu) - V^k(\mu) \geq g(\kappa,\lambda)   \delta_k^2,
	\end{align*}
	which is equivalent to 
	\begin{align}
		\delta_{k+1} \leq  \delta_k - g(\kappa,\lambda)\delta_k^2.
    \label{error_relationship}
	\end{align}
By Assumption \ref{assumption3} and Lemma \ref{lemLambda}, there exist a state $s \in \widetilde{S}$ and $a \notin \calA^*_s$ such that $\pi^k(a|s) > 0 \,\, \forall k \ge 0$ which results in $b^k_s >0$ for all $k \ge 0$. Combined with Assumption \ref{assumption2} and Lemma \ref{keyLemma1}, one can conclude $\delta_k >0$ for all $k \ge 0$. Thus \eqref{error_relationship} is equivalent to 
\begin{align*}
\frac{\delta _k-\delta _{k+1}}{\delta _{k}^{2}}\ge g\left( \kappa ,\lambda \right).
\end{align*}
By Lemma \ref{improvedLowerBound}, the sequence $\{\delta_k\}_{k \ge 0}$ is nonincreasing. Thus
\begin{align*}
\frac{1}{\delta _{k+1}}-\frac{1}{\delta _k}=\frac{\delta _k-\delta _{k+1}}{\delta _{k+1}\delta _k}\ge \frac{\delta _k-\delta _{k+1}}{\delta _{k}^{2}}\ge g\left( \kappa ,\lambda \right),
\end{align*}
and
\begin{align*}
\delta _k=\frac{1}{\frac{1}{\delta _k}}=\frac{1}{\frac{1}{\delta _0}+\sum_{j=0}^{k-1}{\left( \frac{1}{\delta _{j+1}}-\frac{1}{\delta _j} \right)}}\le \frac{1}{\sum_{j=0}^{k-1}{g\left( \kappa ,\lambda \right)}}=\frac{g\left( \kappa ,\lambda \right) ^{-1}}{k},
\end{align*}
which completes the proof.
\end{proof}


\subsection{Proof of Theorem \ref{main result linear}}
\kw{The analysis of the linear convergence relies heavily on the contraction of the optimality measure $b^k_s$ (probability of the output policy $\pi^k$ on non-optimal actions), which is also used in \cite{Khodadadian}.} More precisely, we show that the sub-optimal probabilities $b^k_s$ satisfy (see Lemma~\ref{lemma b})
$$b^k_s - b^{k+1}_s \ge C^k\cdot b^k_s$$ 
 for \kw{sufficiently large} $k$, where the constant $C^k$ depends on the policy $\pi^k$. Plugging this result into Lemma \ref{keyLemma1} yields the linear convergence of Hadamard PG after some calculations.  

\kw{It is worth noting that, based on $b_s^k$,  we can give an intuitive explanation on why softmax PG with a constant step size has at most a sublinear convergence rate while softmax NPG and Hadamard PG with a constant step size can have a local linear convergence rate. First, the update rules of these three algorithms in the policy space are given by
\begin{align*}\text{Softmax PG}:&\quad \pi^{k+1}\left( a|s \right) \propto \pi^k\left( a|s \right) \cdot \exp \{ \eta \cdot \frac{d^k(s) }{1-\gamma} \cdot A^{k}( s,a)  \cdot {\pi^k\left( a|s \right) } \}, \\
\text{Softmax NPG}: &\quad\pi^{k+1}\left( a|s \right) \propto \pi^k\left( a|s \right) \cdot \exp \{ \eta \cdot \frac{d^k(s) }{1-\gamma} \cdot A^{k}( s,a)  \}, \\
\text{Hadamard PG}:&\quad\pi^{k+1}\left( a|s \right) \propto \pi^k\left( a|s \right) \cdot (1+2 \eta \cdot \frac{d^k (s) }{1-\gamma} \cdot A^{k}( s,a))^2,  
\end{align*}
where $\eta $ is a constant step size.
Lemma \ref{keyLemma1} shows that  the convergence  of an algorithm is indeed controlled by the convergence  of the  probability of the output policy on non-optimal actions, that is  $b^{k}_s$.   If the global convergence is established, one has $\forall a\notin \mathcal{A}_s^*: \underset{k\rightarrow \infty}{\lim}\pi^k\left( a|s \right) =0$. Thus the policy update ratio in softmax PG will satisfy  $\forall a \notin \mathcal{A}_s^*:\underset{k\rightarrow \infty}{\lim} \frac{\pi^{k+1}(a|s)}{\pi^k\left( a|s \right)} =1 $ due to the existence of $\pi^k(a|s)$ within the exponential, which roughly implies that $\forall s \in \mathcal {S}: \underset{k\rightarrow \infty}{\lim}\frac{b_s^{k+1}}{b_s^{k}}=1$. Therefore the linear convergence (even locally) of softmax PG cannot be established.  However, for Hadamard PG and softmax NPG, the policy update ratio of any non-optimal action will not be significantly slowed down by the factor $\pi^k\left( a|s \right) $ and is only determined by the advantage function $A^{k}\left( s,a \right)$. Owing to this, we could have $
\underset{k\rightarrow \infty}{\lim} \frac{b_s^{k+1}}{b_s^k}  <1$ since $A^{k}\left( s,a \right)$ is negative for non-optimal actions, and thus the linear convergence can be established.

}
\begin{proof}[Proof of Theorem \ref{main result linear}]
	Applying Lemma \ref{keyLemma1} yields that
	\begin{align*}
		V^\ast(\mu) - V^k(\mu) &\leq \frac{1}{(1-\gamma)^2} \cdot \E[s\sim d_{\mu}^k]{b_s^k},\numberthis\label{eqMain1}
	\end{align*}
    and
    $$\forall s \in \calS : b^k_s \le \frac{\left\|V^*-V^k\right\|_\infty}{\Delta}.$$
    Notice that
    \begin{align}
    \left\|V^*-V^k\right\|_\infty\le \frac{\sum_{s\in \mathcal{S}}{\mu \left( s \right) \left( V^*\left( s \right) -V^k\left( s \right) \right)}}{\tilde{\mu}}=\frac{V^*\left( \mu \right) -V^k\left( \mu \right)}{\tilde{\mu}}\le \frac{g\left( \kappa ,\lambda \right) ^{-1}}{k\tilde{\mu}}\label{max_error_bound}.
    \end{align}
     Thus if 
    \begin{align*}
        \frac{g\left( \kappa ,\lambda \right) ^{-1}}{k\tilde{\mu}}\le \Delta \cdot \min\{\frac{3\delta\kappa}{20\gamma},\frac{1}{2}\}.
    \end{align*}
    or equivalently,
$$k \geq k_0(\delta,\kappa,\gamma) := \frac{g(\kappa,\lambda)^{-1}}{\tilde{\mu} \Delta}\cdot\min\{\frac{3\delta\kappa}{20\gamma},\frac{1}{2}\}^{-1},
$$
one has
\begin{align}
b^k_s \le \frac{1}{2}   \quad \text{and} \quad  \left\|V^*-V^k\right\|_\infty\le \frac{3\kappa\Delta\delta}{20\gamma}.
\label{g1}
\end{align}
Plugging  \eqref{g1} into Lemma \ref{lemma b} yields that
	\begin{align*}
        b_{s}^{k}-b_{s}^{k+1}&\ge \frac{\eta d_{\mu}^{k}(s)}{2(1-\gamma)}\cdot \frac{\min\mathrm{(}1-b_{s}^{k},b_{s}^{k})}{1+\eta ^2\left\| g_{s}^{k} \right\| _{2}^{2}}\cdot \left(3\kappa\Delta -20\gamma \left\|V^*-V^k\right\|_\infty \right),
\\    &\ge \frac{\eta \tilde{\mu}}{2}\cdot \frac{3\kappa\Delta(1-\delta)}{1+\eta ^2\left\| g_{s}^{k} \right\| _{2}^{2}}\cdot b^k_s ,
\\    &\overset{\left( a \right)}{\ge} \frac{\eta \tilde{\mu}}{2}\cdot \frac{3\kappa\Delta(1-\delta)}{2}\cdot b^k_s ,
\\   &=  \frac{3\kappa^2 (1-\gamma)^2 \tilde{\mu}}{16}\cdot \Delta(1-\delta)\cdot b^k_s ,
\\
     &:= \rho(\delta,\kappa,\gamma) \cdot b^k_s,
\numberthis
\label{b_rho}
\end{align*}
where step (a) is due to \eqref{eq 4} and $\kappa \le 2$. By \eqref{b_rho},we have
\begin{align}
b_{s}^{k}\le (1-\rho (\delta,\kappa ,\gamma ))^{k-k_0(\delta,\kappa,\gamma)}b_{s}^{k_0(\kappa,\gamma)} 
\label{main_result_for_b_local},
\end{align}
for all $k\geq k_0(\kappa,\gamma)$ and $s\in \calS$. For $0\leq k\leq k_0(\kappa,\gamma)$, it can be seen that 
	\begin{align*}
		(1-\rho(\kappa,\gamma))^{k-k_0(\kappa,\gamma)} \geq 1 > \underset{s\in \mathcal{S}}{\max}\,\,b_{s}^{k_0(\kappa,\gamma)}.
	\end{align*}
	Thus we have
	\begin{align}
		\label{eqb}
b_{s}^{k}\le (1-\rho (\kappa ,\gamma ))^{k-k_0(\kappa,\gamma)}
	\end{align}
 for all $k \ge 1$. The proof is completed by combining \eqref{eqb} and \eqref{eqMain1}.
 \end{proof}
\section{Proofs of technical lemmas}
\label{proofLemma}

\subsection{Proof of Lemma \ref{keyLemma1}}
\kw{The proof of the upper bound can be found in \cite{Khodadadian}, while the proof of lower bound can be found in \cite{LiuLiWei2023PPG}. Here we provide a unified proof for completeness.}
First a direct calculation yields that 
	\begin{align*}
V^{\ast}(\rho )-V^{\pi}(\rho )&=-\left( V^{\pi}(\rho )-V^{\ast}(\rho ) \right) 
\\
&\overset{(a)}{=}-\frac{1}{1-\gamma}\mathbb{E} _{s\sim d_{\rho}^{\pi}}\mathbb{E} _{a\sim \pi \left( \cdot |s \right)}\left[ A^{\ast}(s,a) \right] 
\\
&\overset{(b)}{=}-\frac{1}{1-\gamma}\mathbb{E} _{s\sim d_{\rho}^{\pi}}\left\{ \sum_{a\notin \calA _{s}^{\ast}}{\pi}(a|s)A^{\ast}(s,a) \right\}
\\
&\overset{(c)}{=}\frac{1}{1-\gamma}\mathbb{E} _{s\sim d_{\rho}^{\pi}}\left\{ \sum_{a\notin \calA _{s}^{\ast}}{\pi}(a|s)\left|A^{\ast}(s,a)\right| \right\},
\end{align*}
	where step (a) follows from Lemma \ref{performanceDifferenceLemma}, step (b) is due to the fact that $A^\ast(s,a)=0$ for any $a\in\calA_s^\ast$, (c) is due to $A^\ast(s,a) \le 0$ for any $a \notin \calA_s^\ast$. Notice that for any $s \in \tilde{S}$ and $a \not\in \calA^\ast$, we have
$$
    \Delta \le \left|A^*(s,a)\right| \le \frac{1}{1-\gamma}.
$$
Thus,
$$
\Delta \cdot \mathbb{E} _{s\sim \rho}\{b^\pi_s\} \le 
\frac{\Delta}{1-\gamma} \mathbb{E} _{s\sim d^\pi_\rho}\{b^\pi_s\} \le
V^\ast(\rho)-V^\pi(\rho) \le 
\frac{1}{(1-\gamma)^2}\mathbb{E} _{s\sim d^\pi_\rho}\{b^\pi_s\}.
$$

\subsection{Proof of Lemma \ref{improvedLowerBound}}
\label{proofImproved}
Notice that
\begin{align*}
		\pi^{k+1}(a|s) &= (\theta_{s,a}^{k+1} )^2\\
		&=\frac{\left( \theta_{s,a}^k + \eta g_{s,a}^k \right)^2}{\twonorm{\vtheta_s^k + \eta \vg_s^k}^2}\\
		&=\frac{\left( \theta_{s,a}^k + \eta g_{s,a}^k \right)^2}{\twonorm{\vtheta_s^k}^2 + \eta^2 \twonorm{\vg_s^k}^2}\\
		&=\theta_{s,a}^k  \frac{\left( 1+ \eta g_{s,a}^k /\theta_{s,a}^k \right)^2}{\twonorm{\vtheta_s^k}^2 + \eta^2 \twonorm{\vg_s^k}^2}\\
		&=\pi^k(a|s)\cdot \frac{\left( 1+ \eta g_{s,a}^k /\theta_{s,a}^k \right)^2}{ 1+ \eta^2 \twonorm{\vg_s^k}^2}, \numberthis\label{eq 6}
	\end{align*}
	where the third line is due to $\la \vtheta_s^k, \vg_s^k\ra= 0$ and the last line is due to $\twonorm{\vtheta_s^k}=1$.  
	Moreover, 
	\begin{align*}
\sum_a{\pi ^{k+1}}(a|s)A^k(s,a)=\frac{\sum_a{\pi ^k}(a|s)\cdot \left( 1+\eta g_{s,a}^{k}/\theta _{s,a}^{k} \right) ^2\cdot A^k(s,a)}{1+\eta ^2\left\| g_{s}^{k} \right\| _{2}^{2}}.\numberthis \label{eq 1}
	\end{align*}
A direct computation yields that
	\begin{align*}
\sum_a{\pi ^k}(a|s)\cdot \left( 1+\eta g_{s,a}^{k}/\theta _{s,a}^{k} \right) ^2\cdot A^k(s,a)&=\sum_a{\pi ^k}(a|s)\left( 1+\frac{2\eta}{1-\gamma}d_{\mu}^{k}(s)A^k(s,a) \right) ^2A^k(s,a)
\\
&=\sum_a{\pi ^k}(a|s)\left( 1+\frac{4\eta}{1-\gamma}d_{\mu}^{k}(s)A^k(s,a)+\left( \frac{2\eta}{1-\gamma}d_{\mu}^{k}(s)A^k(s,a) \right) ^2 \right) A^k(s,a)
\\
&=\sum_a{\pi ^k}(a|s)\left( \frac{4\eta}{1-\gamma}d_{\mu}^{k}(s)A^k(s,a)+\left( \frac{2\eta}{1-\gamma}d_{\mu}^{\pi _{\theta}}(s)A^k(s,a) \right) ^2 \right) A^k(s,a)
\\
&=\sum_a{\pi ^k}(a|s)\left( \frac{4\eta}{1-\gamma}d_{\mu}^{k}(s)+\left( \frac{2\eta}{1-\gamma}d_{\mu}^{k}(s) \right) ^2A^k(s,a) \right) (A^k(s,a))^2
\\
&=\frac{4\eta}{1-\gamma}d_{\mu}^{k}(s)\sum_a{\pi ^k}(a|s)\left( 1+\frac{\eta d_{\mu}^{k}(s)}{1-\gamma}A^k(s,a) \right) (A^k(s,a))^2
\\
&\overset{\left( a \right)}{\ge}\frac{4\eta}{1-\gamma}d_{\mu}^{k}(s)\sum_a{\pi ^k}(a|s)\frac{3}{4}(A^k(s,a))^2
\\
&=\frac{2\eta d_{\mu}^{k}(s)}{1-\gamma}\sum_a{\pi ^k}(a|s)(A^k(s,a))^2,\numberthis \label{eq 2}
	\end{align*}
	where step (a) is due to the following fact, $\eta=\frac{(1-\gamma)^2 \kappa}{4}, \kappa \in(0,2)$ and 
	\begin{align*}
1+\frac{\eta d_{\mu}^{k}(s)}{1-\gamma}A^k(s,a)&\ge 1-\frac{\eta d_{\mu}^{k}(s)}{1-\gamma}\cdot \frac{1}{1-\gamma}
\\
&\ge 1-\frac{\eta}{(1-\gamma )^2}
\\
&=1-\frac{1}{(1-\gamma )^2}\frac{(1-\gamma )^2\kappa}{4}
\\
&=1-\frac{\kappa}{4}
\\
&>\frac{1}{2}.
	\end{align*}
	Plugging \eqref{eq 2} into \eqref{eq 1} yields that
	\begin{align*}
\sum_a{\pi ^{k+1}}(a|s)A^k(s,a)\ge \frac{2\eta d_{\mu}^{k}(s)}{1-\gamma}\cdot \frac{\sum_a{\pi ^k}(a|s)(A^k(s,a))^2}{1+\eta ^2\left\| g_{s}^{k} \right\| _{2}^{2}}. \numberthis \label{eq 3}
	\end{align*}
	Moreover, it can be seen that
	\begin{align*}
1+\eta ^2\left\| g_{s}^{k} \right\| _{2}^{2}&=1+\eta ^2\sum_a{(}g_{s,a}^{k})^2
\\
&=1+\eta ^2\sum_a{\left( \frac{2}{1-\gamma}d_{\mu}^{k}(s)\cdot \theta_{s,a}^k\cdot A^k(s,a) \right) ^2}
\\
&=1+\eta ^2\cdot \left( \frac{2}{1-\gamma}d_{\mu}^{k}(s) \right) ^2\sum_a{\left( \theta _{s,a}^{k} \right) ^2}\left( A^k(s,a) \right) ^2
\\
&\le 1+\eta ^2\cdot \left( \frac{2}{1-\gamma}d_{\mu}^{k}(s) \right) ^2\sum_a{\pi ^k}\left( a|s \right) \frac{1}{(1-\gamma )^2}
\\
&=1+\eta ^2\cdot \left( \frac{2}{1-\gamma}d_{\mu}^{k}(s) \right) ^2\frac{1}{(1-\gamma )^2}
\\
&\le 1+\eta ^2\cdot \left( \frac{2}{1-\gamma} \right) ^2\frac{1}{(1-\gamma )^2}
\\
&=1+\frac{(1-\gamma )^4\kappa ^2}{16}\cdot \frac{4}{(1-\gamma )^4}
\\
&=1+\frac{\kappa ^2}{4}. \numberthis \label{eq 4}
	\end{align*}
	Plugging \eqref{eq 4} into \eqref{eq 3} yields that
	\begin{align*}
\sum_a{\pi ^{k+1}}(a|s)A^k(s,a)&\ge \frac{2\eta d_{\mu}^{k}(s)}{\left( 1-\gamma \right) \left( 1+\frac{\kappa ^2}{4} \right)}\cdot \sum_a{\pi ^k}(a|s)(A^k(s,a))^2
\\
\quad                           &\overset{\left( a \right)}{\ge}\frac{2\eta \tilde{\mu}}{1+\frac{\kappa ^2}{4}}\cdot \sum_a{\pi ^k}(a|s)(A^k(s,a))^2
\\
\quad                           &\overset{\left( b \right)}{=}\frac{2\kappa\tilde{\mu}\left( 1-\gamma \right) ^2}{4+\kappa^2}\cdot \sum_a{\pi ^k}(a|s)(A^k(s,a))^2, \numberthis \label{SingleStateImprove}
	\end{align*}
	where step (a) is due to the fact that
	\begin{align}
		\label{eqD}
		d^{k+1}_{\mu}(s) \geq (1-\gamma) \mu(s)\geq (1-\gamma)\min_s \mu(s)=(1-\gamma)\tilde{\mu},
	\end{align}
        and step (b) is due to  $
\eta =\frac{\left( 1-\gamma \right) ^2\kappa}{4}$. Moreover, by Lemma \ref{performanceDifferenceLemma}, for any  $\rho$,
\begin{align*}
V^{k+1}(\rho )-V^k(\rho )&=\frac{1}{1-\gamma}\sum_s{
	d_{\rho}^{k+1}\left( s \right) \sum_a{\pi ^{k+1}\left( a|s \right) A^k\left( s,a \right)}.
}  \numberthis \label{improvement_on_rho}
\end{align*}
Plugging (\ref{SingleStateImprove}) into (\ref{improvement_on_rho}) yields that
\begin{align*}
V^{k+1}(\rho )-V^k(\rho )&\ge \frac{1}{1-\gamma}\mathbb{E} _{s\sim d_{\rho}^{k+1}}\left\{ \frac{2\kappa \tilde{\mu}\left( 1-\gamma \right) ^2}{4+\kappa^2}\cdot \sum_a{\pi ^k}(a|s)(A^k(s,a))^2 \right\} 
\\
\quad                      &=\frac{2\kappa \tilde{\mu}\left( 1-\gamma \right)}{4+\kappa^2}\sum_s{d_{\rho}^{k+1}\left( s \right) \sum_a{\pi ^k}(a|s)(A^k(s,a))^2}
\\
\quad                      &\ge \frac{2\kappa \tilde{\mu}\left( 1-\gamma \right)}{4+\kappa^2}\sum_{s\in \mathrm{supp}\left( \rho \right)}{d_{\rho}^{k+1}\left( s \right) \sum_a{\pi ^k}(a|s)(A^k(s,a))^2}
\\
\quad                      &\ge \frac{2\kappa \tilde{\mu}\left( 1-\gamma \right) ^2}{4+\kappa^2}\underset{s\in \mathrm{supp}\left( \rho \right)}{\min}\left\{ \tilde{\rho}\left( s \right) \right\} \sum_{s\in \mathrm{supp}\left( \rho \right)}{\sum_a{\pi ^k}(a|s)(A^k(s,a))^2},
\end{align*}
which completes the proof.

\subsection{Proof of Lemma \ref{policyUpdate}}
\label{proofPolicyUpdate}
A straightforward calculation yields that
\begin{align*}
	\pi^{k+1}(a|s) - \pi^k(a|s) &\stackrel{(a)}{=} \pi^k(a|s)\cdot \left( \frac{\left( 1+ \eta g_{s,a}^k /\theta_{s,a}^k \right)^2}{ 1+ \eta^2 \twonorm{\vg_s^k}^2} -1 \right)\\
	&= \pi^k(a|s)\cdot \left( \frac{\left( 1+ \eta g_{s,a}^k /\theta_{s,a}^k \right)^2 - 1 - \eta^2 \twonorm{\vg_s^k}^2}{ 1+ \eta^2 \twonorm{\vg_s^k}^2} \right)\\
	&\stackrel{(b)}{=}\frac{ \pi^k(a|s) }{1+ \eta^2 \twonorm{\vg_s^k}^2} \cdot \frac{4\eta d_{\mu}^k(s)}{1-\gamma} \left( A^k(s,a) + \frac{\eta d_{\mu}^k(s)}{1-\gamma} \left(  (A^{k}(s,a))^2 - \E[a'\sim \pi^k(\cdot|s)]{(A^{k}(s,a'))^2}\right)\right),
\end{align*}
where step (a) is due to \eqref{eq 6} and step (b) is due to 
\begin{align*}
	\left( 1+ \eta g_{s,a}^k /\theta_{s,a}^k \right)^2 - 1 - \eta^2 \twonorm{\vg_s^k}^2 &= 2\eta \frac{ g_{s,a}^k }{\theta_{s,a}^k }  + \eta^2 \left( \frac{ g_{s,a}^k }{\theta_{s,a}^k }\right)^2 - \eta^2 \sum_{a'\in\calA} (g_{s,a'}^k)^2\\
	&=2\eta \cdot  \frac{2d_{\mu}^{k}(s) \cdot  A^{k}(s,a) }{1-\gamma} + \eta^2 \left( \frac{2d_{\mu}^{k}(s) \cdot  A^{k}(s,a) }{1-\gamma}\right)^2 - \eta^2 \sum_{a'\in\calA} (\theta_{s,a'}^k)^2\left( \frac{2d_{\mu}^{k}(s) \cdot  A^{k}(s,a') }{1-\gamma}\right)^2\\
	&= \frac{4\eta d_{\mu}^k(s)}{1-\gamma} \left( A^k(s,a) + \frac{\eta d_{\mu}^k(s)}{1-\gamma} \left(  (A^{k}(s,a))^2 - \E[a'\sim \pi^k(\cdot|s)]{(A^{k}(s,a'))^2}\right)\right). 
\end{align*}

\subsection{Proof of Lemma \ref{lemma b}}
\label{proofLemmab}
Since if $s\notin\widetilde{S}$ (i.e. there does not exist non-optimal actions at $s$), $b_s^k=b_s^{k+1}=0$ and the lemma holds trivially, we will restrict our attention to $s\in\widetilde{S}$. In this case, by the first claim of Lemma~\ref{lemLambda}  (whose proof does not rely on Lemma \ref{lemma b}), we always have $0<b_s^k<1$. To prepare for the proof, we define the two auxiliary policies associated with $\pi^k$ as follows:
    \begin{align*}
\begin{cases}
	\psi^k\left( \cdot |s \right) =\frac{\pi^k \left( \cdot |s \right)}{1-b_{s}^{k}}\mathbf{1}\left[ \cdot \in \mathcal{A} _{s}^{*} \right],\\
	\xi^k\left( \cdot |s \right) =\frac{\pi^k \left( \cdot |s \right)}{b_{s}^{k}}\mathbf{1}\left[ \cdot \notin \mathcal{A} _{s}^{*} \right].\\
\end{cases}
    \end{align*}
    It follows immediately that $\pi^k(\cdot|s)=(1-b_s^k)\psi^k(\cdot|s)+b_s^k\xi^k(\cdot|s)$.
    
By the definition of $b_s^k$, we have 
\begin{align*}
	b_s^k - b_s^{k+1} &=\sum_{a\notin\calA_s^\ast} \left( \pi^k(a|s) - \pi^{k+1}(a|s) \right)\\
	&=\sum_{a\in\calA_s^\ast} \left( \pi^{k+1}(a|s) - \pi^k(a|s)  \right).
\end{align*}
By Lemma \ref{policyUpdate}, one has
\begin{align}
\label{diff pi}
	\pi^{k+1}(a|s) - \pi^k(a|s) &=\frac{ \pi^k(a|s) }{1+ \eta^2 \twonorm{\vg_s^k}^2} \cdot \frac{4\eta d_{\mu}^k(s)}{1-\gamma} \left( A^k(s,a) + \frac{\eta d_{\mu}^k(s)}{1-\gamma} \left(  (A^{k}(s,a))^2 - \E[a'\sim \pi^k(\cdot|s)]{(A^{k}(s,a'))^2}\right)\right).
\end{align}
Below we will derive the expressions for $A^k(s,a)$ and $(A^{k}(s,a))^2 - \E[a'\sim \pi^k(\cdot|s)]{(A^{k}(s,a'))^2}$ based on $b^k_s$. 
\begin{itemize}
	\item Derivation of  $A^k(s,a)$. 
 
 Let $\delta_s^k = V^\ast(s) - V^k(s)$ for any $s\in\calS$. Then one has
\begin{align*}
	Q^k(s,a) &=\E[s'\sim P(\cdot|s,a)]{ r(s,a) + \gamma V^k(s')} \\
	&=\E[s'\sim P(\cdot|s,a)]{ r(s,a) + \gamma V^\ast(s') + \gamma\left( V^k(s') - V^\ast(s') \right)} \\
	&=\E[s'\sim P(\cdot|s,a)]{ r(s,a) + \gamma V^\ast(s') - \gamma \delta_{s'}^k} \\
	&=Q^\ast(s,a) - \gamma  \E[s'\sim P(\cdot|s,a)]{ \delta_{s'}^k} ,\numberthis\label{eq Q}\\
    V^k(s)&=\mathbb{E} _{a'\sim \pi ^k\left( \cdot |s \right)}\left\{ Q^k(s,a') \right\} 
\\
\quad       &=\mathbb{E} _{a'\sim \pi ^k\left( \cdot |s \right)}\left\{ Q^{\ast}(s,a') \right\} -\gamma \mathbb{E} _{a'\sim \pi ^k\left( \cdot |s \right) ,s'\sim P\left( \cdot |s,a' \right)}\left\{ \delta _{s'}^{k} \right\} 
\\
\quad       &\overset{\left( a \right)}{=}\left( 1-b_{s}^{k} \right) \mathbb{E} _{a'\sim \psi ^k\left( \cdot |s \right)}\left\{ Q^{\ast}\left( s,a' \right) \right\} +b_{s}^{k}\mathbb{E} _{a'\sim \xi ^k\left( \cdot |s \right)}\left\{ Q^{\ast}\left( s,a' \right) \right\} -\gamma \mathbb{E} _{a'\sim \pi ^k\left( \cdot |s \right) ,s'\sim P\left( \cdot |s,a' \right)}\left\{ \delta _{s'}^{k} \right\} 
\\
\quad       &\overset{\left( b \right)}{=}\left( 1-b_{s}^{k} \right) V^*\left( s \right) +b_{s}^{k}\mathbb{E} _{a'\sim \xi ^k\left( \cdot |s \right)}\left\{ Q^{\ast}\left( s,a' \right) \right\} -\gamma \mathbb{E} _{a'\sim \pi ^k\left( \cdot |s \right) ,s'\sim P\left( \cdot |s,a' \right)}\left\{ \delta _{s'}^{k} \right\} 
\\
\quad       :&\overset{\left( c \right)}{=}\left( 1-b_{s}^{k} \right) V^*\left( s \right) +b_{s}^{k}\cdot e_{s}^{k}-\gamma \mathbb{E} _{a'\sim \pi ^k\left( \cdot |s \right) ,s'\sim P\left( \cdot |s,a' \right)}\left\{ \delta _{s'}^{k} \right\} 
\\
\quad       &=V^*\left( s \right) - b_{s}^{k}\left( V^*\left( s \right) -e_{s}^{k} \right) -\gamma \mathbb{E} _{a'\sim \pi ^k\left( \cdot |s \right) ,s'\sim P\left( \cdot |s,a' \right)}\left\{ \delta _{s'}^{k} \right\}, \numberthis\label{eq V}
\end{align*}
where we have used the definition $b_s^k = \sum_{a\notin \calA_s^\ast} \pi^k(a|s)$, step (a) follows from the definitions of $\psi^k$ and $\xi^k$, step (b) holds since $Q^\ast(s,a') = V^\ast(s)$ for any $a\in\calA_s^\ast$ and $\text{supp}(\psi^k(\cdot|s)) \in \calA^*_s$, and we define $e_s^k = \E[a' \sim \xi^k(\cdot|s)]{Q^\ast(s,a')}$ in  step (c). Therefore, we have
\begin{align*}
	A^{k}(s,a) &= Q^k(s,a) - V^k(s)\\
	&=Q^\ast(s,a) - \gamma  \E[s'\sim P(\cdot|s,a)]{ \delta_{s'}^k}  - \left(V^\ast(s) - b_s^k (V^\ast(s) - e_s^k) - \gamma  \E[s'\sim P(\cdot|s,a'), a'\sim\pi^k(\cdot|s)]{ \delta_{s'}^k}  \right)\\
	&=A^\ast(s,a) + b_s^k (V^\ast(s) - e_s^k) +\underbrace{ \gamma  \left( \E[a'\sim\pi^k(\cdot|s),s'\sim P(\cdot|s,a')]{ \delta_{s'}^k}  - \E[s'\sim P(\cdot|s,a)]{ \delta_{s'}^k}   \right) }_{\ell_{s,a}^k}\\
	:&=A^\ast(s,a) + b_s^k (V^\ast(s) - e_s^k) + \ell_{s,a}^k.
\end{align*}
Furthermore, for any $a\in\calA_s^\ast$, one has
\begin{align}
	\label{eq 10}
	A^{k}(s,a) =b_s^k (V^\ast(s) - e_s^k) + \ell_{s,a}^k.
\end{align}

By the definition of $\ell_{s,a}^k$, it can be seen that
\begin{align*}
	\E[a\sim \pi^k(\cdot|s)]{\ell_{s,a}^k} &=\gamma   \E[a\sim \pi^k(\cdot|s)]{ \left( \E[s'\sim P(\cdot|s,a'), a'\sim\pi^k(\cdot|s)]{ \delta_{s'}^k}  - \E[s'\sim P(\cdot|s,a)]{ \delta_{s'}^k}   \right)  }\\
	&=\gamma   \E[a\sim \pi^k(\cdot|s)]{ \left( \E[s'\sim P(\cdot|s)]{ \delta_{s'}^k}  - \E[s'\sim P(\cdot|s,a)]{ \delta_{s'}^k}   \right)  }\\
	&=0.
\end{align*}

\item Derivation of $(A^{k}(s,a))^2 - \E[a'\sim \pi^k(\cdot|s)]{(A^{k}(s,a'))^2}$.

A direct computation yields that 
\begin{align*}
	&(A^k(s,a))^2 - \E[a'\sim \pi^k(\cdot|s)]{(A^k(s,a'))^2 }\\
	&= \left( Q^k(s,a) - V^k(s) \right)^2 - \E[a'\sim \pi^k(\cdot|s)]{\left( Q^k(s,a') - V^k(s) \right)^2}\\
	&= \left( Q^k(s,a)  \right)^2 - 2 Q^k(s,a) \cdot V^k(s) + \left(V^k(s)\right)^2 -  \E[a'\sim \pi^k(\cdot|s)]{\left( Q^k(s,a')  \right)^2 - 2 Q^k(s,a') \cdot V^k(s) + \left(V^k(s)\right)^2}\\
	&= \left( Q^k(s,a)  \right)^2 - 2 Q^k(s,a) \cdot V^k(s) + \left(V^k(s)\right)^2 -  \E[a'\sim \pi^k(\cdot|s)]{\left( Q^k(s,a')  \right)^2 }\\
 &\quad + 2 \E[a'\sim \pi^k(\cdot|s)]{Q^k(s,a') }\cdot V^k(s)- \left(V^k(s)\right)^2\\
&= \left( Q^k(s,a)  \right)^2 -\E[a'\sim \pi^k(\cdot|s)]{\left( Q^k(s,a')  \right)^2 }- 2 Q^k(s,a) \cdot V^k(s) + 2 (V^k(s))^2\\
&= \left( Q^k(s,a)  \right)^2 -\E[a'\sim \pi^k(\cdot|s)]{\left( Q^k(s,a')  \right)^2 }- 2 \left( Q^k(s,a)- V^k(s)\right) \cdot V^k(s) \\
&= \left( Q^k(s,a)  \right)^2 -\E[a'\sim \pi^k(\cdot|s)]{\left( Q^k(s,a')  \right)^2 }- 2A^k(s,a)\cdot V^k(s). \numberthis\label{eq 7}
\end{align*}
For the first term in \eqref{eq 7}, using \eqref{eq Q} yields that
\begin{align*}
	 &\left( Q^k(s,a)  \right)^2 -\E[a'\sim \pi^k(\cdot|s)]{\left( Q^k(s,a')  \right)^2 } \\
	 &= \left(Q^\ast(s,a) - \gamma  \E[s'\sim P(\cdot|s,a)]{ \delta_{s'}^k} \right)^2 - \E[a'\sim \pi^k(\cdot|s)]{\left(Q^\ast(s,a') - \gamma  \E[s'\sim P(\cdot|s,a)]{ \delta_{s'}^k}\right)^2}\\
	 &=\left(Q^\ast(s,a) \right)^2 - \E[a'\sim \pi^k(\cdot|s)]{\left(Q^\ast(s,a')  \right)^2}\\
	 &\quad + \gamma^2 \left( \E[s'\sim P(\cdot|s,a)]{ \delta_{s'}^k} \right)^2 -  \gamma^2\E[a'\sim \pi^k(\cdot|s)]{\left( \E[s'\sim P(\cdot|s,a)]{ \delta_{s'}^k}\right)^2}\\
	 &\quad + 2\gamma \E[a'\sim \pi^k(\cdot|s)]{ Q^\ast(s,a') \cdot \E[s'\sim P(\cdot|s,a)]{ \delta_{s'}^k} } - 2\gamma Q^\ast(s,a) \cdot \E[s'\sim P(\cdot|s,a)]{ \delta_{s'}^k}\\
	 :&=\left(Q^\ast(s,a) \right)^2 - \E[a'\sim \pi^k(\cdot|s)]{\left(Q^\ast(s,a')  \right)^2} + \tilde{\ell}_{s,a}^k. \numberthis \label{eq 8}
\end{align*}
By the definitions of $\psi^k$ and $\xi^k$, one has
\begin{align*}
\mathbb{E} _{a'\sim \pi ^k\left( \cdot |s \right)}\left\{ \left( Q^{\ast}(s,a') \right) ^2 \right\} &=\left( 1-b_{s}^{k} \right) \mathbb{E} _{a'\sim \psi ^k\left( \cdot |s \right)}\left\{ \left( Q^{\ast}(s,a') \right) ^2 \right\} +b_{s}^{k} \cdot \mathbb{E} _{a'\sim \xi ^k\left( \cdot |s \right)}\left\{ \left( Q^{\ast}(s,a') \right) ^2 \right\} 
\\
\quad                             &\overset{(a)}{:=}\left( 1-b_{s}^{k} \right) \mathbb{E} _{a'\sim \psi ^k\left( \cdot |s \right)}\left\{ \left( V^*\left( s \right) \right) ^2 \right\} +b_{s}^{k}\cdot \left( d_{s}^{k} \right) ^2
\\
\quad                              &=\left( 1-b_{s}^{k} \right) \left( V^*\left( s \right) \right) ^2+b_{s}^{k}\cdot \left( d_{s}^{k} \right)^2,\numberthis \label{eq 9}
\end{align*}
where we define $
d_{s}^{k}:=\sqrt{\mathbb{E} _{a'\sim \xi ^k\left( \cdot |s \right)}\left\{ \left( Q^{\ast}(s,a') \right) ^2 \right\}}$ in step (a).
Plugging \eqref{eq 8} and \eqref{eq 9} into \eqref{eq 7} yields that for any $(s,a)$,
\begin{align*}
	&(A^k(s,a))^2 - \E[a'\sim \pi^k(\cdot|s)]{(A^k(s,a'))^2 }\\
 &= \left( Q^k(s,a)  \right)^2 -\E[a'\sim \pi^k(\cdot|s)]{\left( Q^k(s,a')  \right)^2 }- 2A^k(s,a)\cdot V^k(s)\\
	&= \left(Q^\ast(s,a) \right)^2 - \E[a'\sim \pi^k(\cdot|s)]{\left(Q^\ast(s,a')  \right)^2} + \tilde{\ell}_{s,a}^k - 2A^k(s,a)\cdot V^k(s)\\
	&=\left(Q^\ast(s,a) \right)^2 - (1-b_s^k) \left(V^\ast(s)  \right)^2 - b_s^k\cdot (d_s^k)^2 + \tilde{\ell}_{s,a}^k - 2A^k(s,a)\cdot V^k(s)\\
	&= \left(Q^\ast(s,a) \right)^2 -  \left(V^\ast(s)  \right)^2 + b_s^k\cdot\left(  \left(V^\ast(s)  \right)^2  - (d_s^k)^2  \right)+ \tilde{\ell}_{s,a}^k - 2A^k(s,a)\cdot V^k(s).
\end{align*}
Moreover, for any $a\in\calA_s^\ast$, one has
\begin{align}
	\label{eq 11}
	&(A^k(s,a))^2 - \E[a'\sim \pi^k(\cdot|s)]{(A^k(s,a'))^2 } \notag\\
	&=b_s^k\cdot\left(  \left(V^\ast(s)  \right)^2  - (d_s^k)^2  \right)+ \tilde{\ell}_{s,a}^k - 2A^k(s,a)\cdot V^k(s)\notag\\
	&= b_s^k\cdot\left(  \left(V^\ast(s)  \right)^2  - (d_s^k)^2  \right)+ \tilde{\ell}_{s,a}^k - 2\left(  b_s^k (V^\ast(s) - e_s^k) + \ell_{s,a}^k\right) \cdot V^k(s)\notag\\
	&=b_s^k\cdot\left(  \left(V^\ast(s)  \right)^2  - (d_s^k)^2 - 2  (V^\ast(s) - e_s^k) \cdot V^k(s) \right)+ \tilde{\ell}_{s,a}^k - 2 \ell_{s,a}^k \cdot V^k(s).
\end{align}
\end{itemize}
Plugging \eqref{eq 10} and \eqref{eq 11} into \eqref{diff pi} yields that
\begin{align*}
	&\pi^{k+1}(a|s) - \pi^k(a|s) \\
 &=\frac{ \pi^k(a|s) }{1+ \eta^2 \twonorm{\vg_s^k}^2} \cdot \frac{4\eta d_{\mu}^k(s)}{1-\gamma} \left( A^k(s,a) + \frac{\eta d_{\mu}^k(s)}{1-\gamma} \left(  (A^{k}(s,a))^2 - \E[a'\sim \pi^k(\cdot|s)]{(A^{k}(s,a'))^2}\right)\right)\\
	&= \frac{ \pi^k(a|s) }{1+ \eta^2 \twonorm{\vg_s^k}^2} \cdot \frac{4\eta d_{\mu}^k(s)}{1-\gamma} \left( b_s^k (V^\ast(s) - e_s^k) + \ell_{s,a}^k \right)\\
	&\quad +   \frac{ \pi^k(a|s) }{1+ \eta^2 \twonorm{\vg_s^k}^2} \cdot\frac{4\eta d_{\mu}^k(s)}{1-\gamma} \cdot \frac{\eta d_{\mu}^k(s)}{1-\gamma} \left( b_s^k\cdot\left(  \left(V^\ast(s)  \right)^2  - (d_s^k)^2 - 2  (V^\ast(s) - e_s^k) \cdot V^k(s) \right)+ \tilde{\ell}_{s,a}^k - 2 \ell_{s,a}^k \cdot V^k(s)\right) \\
	&= \frac{ \pi^k(a|s) }{1+ \eta^2 \twonorm{\vg_s^k}^2} \cdot\frac{4\eta d_{\mu}^k(s)}{1-\gamma}  \cdot  b_s^k (V^\ast(s) - e_s^k)  - \frac{ \pi^k(a|s) }{1+ \eta^2 \twonorm{\vg_s^k}^2} \cdot\frac{4\eta d_{\mu}^k(s)}{1-\gamma} \frac{\eta d_{\mu}^k(s)}{1-\gamma} \cdot b_s^k\cdot 2  (V^\ast(s) - e_s^k) \cdot V^k(s) \\
	&\quad +\frac{ \pi^k(a|s) }{1+ \eta^2 \twonorm{\vg_s^k}^2} \cdot \frac{4\eta d_{\mu}^k(s)}{1-\gamma} \cdot \ell_{s,a}^k +  \frac{ \pi^k(a|s) }{1+ \eta^2 \twonorm{\vg_s^k}^2} \cdot\frac{4\eta d_{\mu}^k(s)}{1-\gamma} \frac{\eta d_{\mu}^k(s)}{1-\gamma} \left( b_s^k\cdot\left(  \left(V^\ast(s)  \right)^2  - (d_s^k)^2  \right)+ \tilde{\ell}_{s,a}^k - 2 \ell_{s,a}^k \cdot V^k(s)\right) \\
	&=\frac{ \pi^k(a|s) }{1+ \eta^2 \twonorm{\vg_s^k}^2} \cdot\frac{4\eta d_{\mu}^k(s)}{1-\gamma}  \cdot  b_s^k \cdot \left[ (V^\ast(s) - e_s^k)  \left( 1 - \frac{2\eta d_{\mu}^k(s)}{1-\gamma} \cdot  V^k(s)\right) + \frac{\eta d_{\mu}^k(s)}{1-\gamma}  \left(  \left(V^\ast(s)  \right)^2  - (d_s^k)^2  \right) \right]\\
	&\quad + \frac{ \pi^k(a|s) }{1+ \eta^2 \twonorm{\vg_s^k}^2} \cdot \frac{4\eta d_{\mu}^k(s)}{1-\gamma} \left[ \ell_{s,a}^k \left(1- \frac{2\eta d_{\mu}^k(s)}{1-\gamma} \cdot  V^k(s) \right) + \frac{\eta d_{\mu}^k(s)}{1-\gamma} \cdot \tilde{\ell}_{s,a}^k  \right]\\
	:&=\frac{ \pi^k(a|s) }{1+ \eta^2 \twonorm{\vg_s^k}^2} \cdot\frac{4\eta d_{\mu}^k(s)}{1-\gamma}  \cdot  \left( b_s^k \cdot z_s^k +  \phi_{s,a}^k \right) 
\end{align*}
for any $a\in\calA_s^\ast$. Thus,
\begin{align*}
	b_s^k - b_s^{k+1} &=\sum_{a\in\calA_s^\ast} \left( \pi^{k+1}(a|s) - \pi^k(a|s)  \right)\\
	&=\frac{4\eta d_{\mu}^k(s)}{1-\gamma}  \cdot  \frac{1}{1+ \eta^2 \twonorm{\vg_s^k}^2} \sum_{a\in\calA_s^\ast} \pi^k(a|s) \left( b_s^k \cdot z_s^k +  \phi_{s,a}^k \right) \\
	&=\frac{4\eta d_{\mu}^k(s)}{1-\gamma}  \cdot  \frac{1}{1+ \eta^2 \twonorm{\vg_s^k}^2} \left( (1-b_s^k)b_s^k \cdot z_s^k + \sum_{a\in\calA_s^\ast} \pi^k(a|s)  \phi_{s,a}^k \right).\numberthis\label{eqB}
\end{align*}
Next we will show
\begin{align*}
	T_s:&= (1-b_s^k)b_s^k \cdot z_s^k + \sum_{a\in\calA_s^\ast} \pi^k(a|s)  \phi_{s,a}^k \\
	&\geq  (1-b_s^k)b_s^k \cdot z_s^k - \left| \sum_{a\in\calA_s^\ast} \pi^k(a|s)  \phi_{s,a}^k \right|.\numberthis\label{eqT}
\end{align*}
Recall the definitions of $\phi_{s,a}^k,  \ell_{s,a}^k , \tilde{\ell}_{s,a}^k$:
\begin{align*}
	\phi_{s,a}^k &= \ell_{s,a}^k \left(1- \frac{2\eta d_{\mu}^k(s)}{1-\gamma} \cdot  V^k(s) \right) + \frac{\eta d_{\mu}^k(s)}{1-\gamma} \cdot \tilde{\ell}_{s,a}^k,\\
	\ell_{s,a}^k &=\gamma  \left( \E[s'\sim P(\cdot|s,a'), a'\sim\pi^k(\cdot|s)]{ \delta_{s'}^k}  - \E[s'\sim P(\cdot|s,a)]{ \delta_{s'}^k}   \right),\\
	\tilde{\ell}_{s,a}^k &= \gamma^2 \left( \E[s'\sim P(\cdot|s,a)]{ \delta_{s'}^k} \right)^2 -  \gamma^2\E[a'\sim \pi^k(\cdot|s)]{\left( \E[s'\sim P(\cdot|s,a)]{ \delta_{s'}^k}\right)^2}\\
	&\quad + 2\gamma \E[a'\sim \pi^k(\cdot|s)]{ Q^\ast(s,a') \cdot \E[s'\sim P(\cdot|s,a)]{ \delta_{s'}^k} } - 2\gamma Q^\ast(s,a) \cdot \E[s'\sim P(\cdot|s,a)]{ \delta_{s'}^k}.
\end{align*}
We have 
\begin{align*}
	\sum_{a\in\calA} \pi^k(a|s) \ell_{s,a}^k  &= 0,\\
	\sum_{a\in\calA} \pi^k(a|s) \tilde{\ell}_{s,a}^k&=0,\\
	\sum_{a\in\calA} \pi^k(a|s)  \phi_{s,a}^k  &=0.
\end{align*}
It follows that
\begin{align*}
	 \left| \sum_{a\in\calA_s^\ast} \pi^k(a|s)  \phi_{s,a}^k \right| &\leq  \sum_{a\in\calA_s^\ast} \pi^k(a|s)   \left| \phi_{s,a}^k \right|\leq \sum_{a\in\calA_s^\ast} \pi^k(a|s)   \cdot \max_{s,a}\left| \phi_{s,a}^k \right| =(1-b_s^k)  \max_{s,a}\left| \phi_{s,a}^k \right| 
\end{align*}
and 
\begin{align*}
	\left| \sum_{a\in\calA_s^\ast} \pi^k(a|s)  \phi_{s,a}^k \right| & =\left| \sum_{a\in\calA}\pi^k(a|s)  \phi_{s,a}^k -  \sum_{a\notin\calA_s^\ast} \pi^k(a|s)  \phi_{s,a}^k \right| \\
	&=\left| \sum_{a\notin\calA_s^\ast} \pi^k(a|s)  \phi_{s,a}^k \right| \\
	&\leq b_s^k \cdot \max_{s,a}\left| \phi_{s,a}^k \right| . \numberthis \label{eqPhi}
\end{align*}
Combining them together yields that
\begin{align}
	\label{eq14}
	\left| \sum_{a\in\calA_s^\ast} \pi^k(a|s)  \phi_{s,a}^k \right|  \leq \min(b_s^k, 1-b_s^k)\cdot \max_{s,a}\left|\phi_{s,a}^k \right| .
\end{align}

On the other hand, one has
\begin{align}
	\label{eq15}
	(1-b_s^k)\cdot b_s^k &=\min(1-b_s^k, b_s^k)\cdot\max (1-b_s^k, b_s^k)\\
	&\geq \frac{1}{2}\min(1-b_s^k, b_s^k),
\end{align}
where the second line is due to the fact that $\max(1-x,x)\geq 1/2$ for any $x\in[0,1]$. Plugging \eqref{eq14} and \eqref{eq15} into \eqref{eqT} yields that
\begin{align*}
	T_s&\geq  (1-b_s^k)b_s^k \cdot z_s^k - \left| \sum_{a\in\calA_s^\ast} \pi^k(a|s)  \phi_{s,a}^k \right|\\
	&\geq \min(1-b_s^k, b_s^k)\left( \frac{z_s^k}{2} -\max_{s,a}\left|\phi_{s,a}^k \right|  \right). \numberthis\label{eqT2}
\end{align*}
In what follows, we will control $z_s^k$ and $|\phi_{s,a}^k |$, respectively. 

\begin{itemize}
	\item Bounding $z_s^k$.

Recall the definition of 
\begin{align*}
	z_s^k &= (V^\ast(s) - e_s^k)  \left( 1 - \frac{2\eta d_{\mu}^k(s)}{1-\gamma} \cdot  V^k(s)\right) + \frac{\eta d_{\mu}^k(s)}{1-\gamma}  \left(  \left(V^\ast(s)  \right)^2  - (d_s^k)^2  \right).
\end{align*}
According to the Jensen inequality, one has
\begin{align*}
	\widehat{V}^*(s) \le e_s^k & = \E[a' \sim \xi^k(\cdot|s)]{Q^*(s,a')}
        \le \sqrt{\mathbb{E} _{a'\sim \xi ^k\left( \cdot |s \right)}\left\{ \left( Q^{\ast}(s,a') \right) ^2 \right\}}=d_{s}^{k} 
        \le  \widetilde{V}^*(s),
\end{align*}
which implies that  $V^\ast(s) - e_s^k \ge V^\ast(s) - d_s^k \ge V^\ast(s) - \widetilde{V}^*(s) > 0$. Thus,
\begin{align*}
	z_s^k &= (V^\ast(s) - e_s^k)  \left( 1 - \frac{2\eta d_{\mu}^k(s)}{1-\gamma} \cdot  V^k(s)\right) + \frac{\eta d_{\mu}^k(s)}{1-\gamma}  \left(  \left(V^\ast(s)  \right)^2  - (d_s^k)^2  \right)\\
 &\ge (V^\ast(s) - d_s^k)  \left( 1 - \frac{2\eta d_{\mu}^k(s)}{1-\gamma} \cdot  V^k(s)\right) + \frac{\eta d_{\mu}^k(s)}{1-\gamma}  \left(  \left(V^\ast(s)  \right)^2  - (d_s^k)^2  \right)\\
  &= (V^\ast(s) - d_s^k)  \left( 1 + \frac{\eta d_{\mu}^k(s)}{1-\gamma}  \left( V^\ast(s) + d^k_s - 2V^k(s)  \right) \right ) \\
  &= (V^\ast(s) - d_s^k)  \left( 1 - \frac{\eta d_{\mu}^k(s)}{1-\gamma}  \left( V^\ast(s) - d^k_s \right) \right ) \\
  &= (V^\ast(s) - d_s^k)  \left( 1 - \frac{\eta d_{\mu}^k(s)}{(1-\gamma)^2} \right ) \\
  &= (V^\ast(s) - d_s^k)  \left( 1 - \frac{\kappa}{4} \right ) \\
&= (V^\ast(s) - \tilde{V}^\ast(s))  \left( \kappa - \frac{\kappa}{4} \right ) \\
    &\ge \frac{3\kappa\Delta}{4} \numberthis \label{eqM}
\end{align*}

\item Bounding $|\phi_{s,a}^k |$. By the definition, we have
	\begin{align*}
		|\phi_{s,a}^k |  &= \left| \ell_{s,a}^k \left(1- \frac{2\eta d_{\mu}^k(s)}{1-\gamma} \cdot  V^k(s) \right) + \frac{\eta d_{\mu}^k(s)}{1-\gamma} \cdot \tilde{\ell}_{s,a}^k \right|\\
		&\leq \left| \ell_{s,a}^k \right| \cdot \left|  1- \frac{2\eta d_{\mu}^k(s)}{1-\gamma} \cdot  V^k(s) \right| + \frac{\eta d_{\mu}^k(s)}{1-\gamma} \cdot \left|\tilde{\ell}_{s,a}^k \right|\\
		&\leq \left| \ell_{s,a}^k \right| + \frac{\eta}{1-\gamma} \cdot \left|\tilde{\ell}_{s,a}^k \right|.
	\end{align*} 
We will control $ \left| \ell_{s,a}^k \right| $ and $\left|\tilde{\ell}_{s,a}^k \right|$, respectively.
\begin{itemize}
	\item Bounding $ \left| \ell_{s,a}^k \right| $. Recall the definition
	\begin{align*}
		\ell_{s,a}^k =\gamma  \left( \E[s'\sim P(\cdot|s,a'), a'\sim\pi^k(\cdot|s)]{ \delta_{s'}^k}  - \E[s'\sim P(\cdot|s,a)]{ \delta_{s'}^k}   \right) .
	\end{align*}
A direct computation yields that
\begin{align*}
	\left| \ell_{s,a}^k \right|  &=\gamma \left| \E[s'\sim P(\cdot|s,a)]{ \delta_{s'}^k}  - \E[a'\sim\pi^k(\cdot|s)]{\E[s'\sim P(\cdot|s,a')]{ \delta_{s'}^k} }\right|\\
	&\leq \gamma \left| \E[s'\sim P(\cdot|s,a)]{ \delta_{s'}^k}   \right|\\
	&\leq \gamma \E[s'\sim P(\cdot|s,a)]{ \left|  \delta_{s'}^k \right|}  \\
	&=\gamma \E[s'\sim P(\cdot|s,a)]{ \left|  V^\ast(s') - V^k(s') \right|}  \\
	&\leq \gamma\left\|V^*-V^k\right\|_\infty.
\end{align*}
	\item Bounding $\left|\tilde{\ell}_{s,a}^k \right|$. Similarly, one has
	\begin{align*}
		\left|\tilde{\ell}_{s,a}^k \right| &\leq \gamma^2 \left|  \left( \E[s'\sim P(\cdot|s,a)]{ \delta_{s'}^k} \right)^2 - \E[a'\sim \pi^k(\cdot|s)]{\left( \E[s'\sim P(\cdot|s,a)]{ \delta_{s'}^k}\right)^2}\right|\\
		&\quad + 2\gamma \left| \E[a'\sim \pi^k(\cdot|s)]{ Q^\ast(s,a') \cdot \E[s'\sim P(\cdot|s,a)]{ \delta_{s'}^k} } - Q^\ast(s,a) \cdot \E[s'\sim P(\cdot|s,a)]{ \delta_{s'}^k}\right|\\
		&\leq \gamma^2 \left\|V^*-V^k\right\|_\infty^2 + 2\gamma V^\ast(s)\cdot \left\|V^*-V^k\right\|_\infty\\
		&\leq \left( \gamma^2 \cdot \left\|V^*-V^k\right\|_\infty +2\gamma \cdot \frac{1}{1-\gamma} \right)\cdot \left\|V^*-V^k\right\|_\infty\\
		&\leq \left( \gamma^2 \cdot \max_s (V^\ast(s) ) +2\gamma \cdot \frac{1}{1-\gamma} \right)\cdot \left\|V^*-V^k\right\|_\infty\\
		&\leq \left( \gamma^2 \cdot \frac{1}{1-\gamma}+2\gamma \cdot \frac{1}{1-\gamma} \right)\cdot \left\|V^*-V^k\right\|_\infty\\
		&\leq \frac{\gamma^2+2\gamma}{1-\gamma}\left\|V^*-V^k\right\|_\infty.
	\end{align*}
\end{itemize}
Combining them together, we have
\begin{align*}
	|\phi_{s,a}^k |  &\leq \left| \ell_{s,a}^k \right| + \frac{\eta}{1-\gamma} \cdot \left|\tilde{\ell}_{s,a}^k \right|\\
	&\leq \left( \gamma + \frac{\eta}{1-\gamma}  \frac{ \gamma^2+2\gamma}{1-\gamma} \right)\left\|V^*-V^k\right\|_\infty\\
	&\leq \left( \gamma + \frac{\eta}{1-\gamma}  \frac{ 3\gamma}{1-\gamma} \right)\left\|V^*-V^k\right\|_\infty\\
	&=\left( \gamma + \frac{(1-\gamma)^2\kappa}{4}\frac{1}{1-\gamma}  \frac{ 3\gamma}{1-\gamma} \right)\left\|V^*-V^k\right\|_\infty\\
	&\leq \frac{5}{2}\gamma \left\|V^*-V^k\right\|_\infty.\numberthis\label{eqM2}
\end{align*} 
\end{itemize}
Plugging \eqref{eqM} and \eqref{eqM2} into \eqref{eqT2} yields that
\begin{align}
	T_s&\geq \min(1-b_s^k, b_s^k)\left( \frac{z_s^k}{2} -\max_{s,a}\left|\phi_{s,a}^k \right|  \right)
    \label{lower_bound_of_t_for_all_s}
 \\
	&\geq \min(1-b_s^k, b_s^k)\left(\frac{3\kappa\Delta}{8} -\frac{5}{2}\gamma \left\|V^*-V^k\right\|_\infty \right).  \label{lower_bound_of_t_for_s_in_tilde_S}
\end{align} 
Submitting it into \eqref{eqB} completes the proof.

\subsection{Proof of Lemma \ref{lemLambda}}
\label{proofLemLambda}
Firstly, we show that $1- b_s^k >0$. 
Recall that
\begin{align*}
	g_{s,a}^k= \frac{2}{1-\gamma} d_{\mu}^{\pi^k}(s)  \theta_{s,a}^k A^{k}(s,a) \text{ and } \vtheta_s^{k+1} = \frac{\vtheta_s^k + \eta \vg_s^k}{\twonorm{\vtheta_s^k + \eta \vg_s^k }}\text{ for }s\in\calS.
\end{align*}
For any state $s\in\calS$, a direct computation yields that
\begin{align*}
	\pi^{k+1}(a|s) &= \pi^k(a|s)\cdot \frac{\left( 1+ \eta g_{s,a}^k /\theta_{s,a}^k \right)^2}{ 1+ \eta^2 \twonorm{\vg_s^k}^2} \\
	&= \frac{\pi^k(a|s)}{ 1+ \eta^2 \twonorm{\vg_s^k}^2} \left( 1+ \eta g_{s,a}^k /\theta_{s,a}^k \right)^2\\
	&\geq \frac{\pi^k(a|s)}{1+\eta^2 \twonorm{\vg_s^k}^2} \left( 1 + \frac{4\eta}{1-\gamma} d_{\mu}^k(s)    A^{k}(s,a) \right)\\
	&\geq \frac{\pi^k(a|s)}{1+\eta^2 \twonorm{\vg_s^k}^2} \left( 1 -\frac{4\eta}{1-\gamma} d_{\mu}^k(s)   \left| A^{k}(s,a)\right| \right)\\
	&\stackrel{(a)}{\geq}\frac{\pi^k(a|s)}{1+\eta^2 \twonorm{\vg_s^k}^2} \left( 1 -\frac{4\eta}{(1-\gamma)^2}  \right)\\
	&\stackrel{(b)}{=}\frac{\pi^k(a|s)}{1+\eta^2 \twonorm{\vg_s^k}^2} \left( 1 -\kappa \right),
\end{align*}
where step (a) follows from $d_{\mu}^k\leq 1$ and $|A^{k}(s,a)|\leq \frac{1}{1-\gamma}$ and step (b) is due to $\eta = \frac{(1-\gamma)^2\kappa}{4}$. Since $\pi^0(a|s)>0$, it can be seen that $\pi^k(a|s)>0$ for all $k\geq 1$. Thus we have
\begin{align*}
	1-b_s^k =1 - \sum_{a\notin \calA_s^\ast} \pi^k(a|s) =  \sum_{a\in \calA_s^\ast} \pi^k(a|s) >0
\end{align*}
for all $k\geq 0$. Moreover, Lemma \ref{lemma b} implies that the sequence $\{b_s^k\}_{k\geq k_0}$ is nonincreasing when $
\underset{s\in \mathcal{S}}{\max}\left( V^*\left( s \right) -V^k\left( s \right) \right)$ is sufficiently small. Thus by Theorem \ref{main result global}, one can conclude that  there exist a constant $k_0$ such that $\{b_s^k\}_{k\geq k_0}$ is nonincreasing for $k\ge k_0$ and thus we have $\lambda := \underset{k\ge0}{\inf} \left\{ \underset{ s\in\calS}{\min} (1-b_s^k) \right\} >0$.


\section{Discussion}
\label{conclusion}
In this paper, we consider  Hadamard parameterization in  policy gradient  and the local  linear convergence of Hadamard PG is established in the tabular setting provided the exact gradient can be accessed.  For future work, \kw{we would like to see whether  the analysis can be extended to a broader range of step sizes.}
In addition, it is interesting to study the convergence of the algorithm in the function approximation setting as well as in random sampling setting. 
\kw{Moreover, it is likely to accelerate the algorithm using the idea of  natural gradient  and meanwhile establish the convergence of the corresponding algorithm. }


\bibliographystyle{plain}
\bibliography{refs}
\end{document}